\documentclass[9pt,shortpaper,twoside,web]{ieeecolor}
\usepackage{generic}
\usepackage{cite}
\usepackage{amsmath,amssymb,amsfonts}
\usepackage{graphicx}
\usepackage{xspace}
\usepackage{textcomp}
\usepackage{mathtools}

\usepackage{booktabs}
\usepackage{multirow}
\usepackage{multicol}
\usepackage{algorithm}
\usepackage{xfrac}
\usepackage{algpseudocode}
\usepackage[linesnumbered,ruled,vlined,algo2e]{algorithm2e}

\DeclareMathOperator*{\argmin}{arg\,min}

\mathtoolsset{showonlyrefs}

\newcommand{\R}{\ensuremath{\mathbb{R}}\xspace}
\newcommand{\X}{\ensuremath{\mathcal{X}}\xspace}
\newcommand{\Z}{\ensuremath{\mathcal{Z}}\xspace}
\newcommand{\V}{\ensuremath{\mathcal{V}}\xspace}
\newcommand{\G}{\ensuremath{\mathcal{G}}\xspace}
\newcommand{\E}{\ensuremath{\mathcal{E}}\xspace}
\newcommand{\col}[2]{\textnormal{col}(#1,#2)}
\newcommand{\tr}{\textnormal{true}}
\newcommand{\ops}[1]{\textnormal{ops}(#1)}

\newtheorem{problem}{Problem}
\newtheorem{assumption}{Assumption}
\newtheorem{lemma}{Lemma}
\newtheorem{definition}{Definition}
\newtheorem{theorem}{Theorem}
\newtheorem{corollary}{Corollary}

\def\BibTeX{{\rm B\kern-.05em{\sc i\kern-.025em b}\kern-.08em
    T\kern-.1667em\lower.7ex\hbox{E}\kern-.125emX}}
\begin{document}

\title{A Distributed Asynchronous Generalized Momentum \\Algorithm Without Delay Bounds}
\author{Ellie Pond$^{1*}$, Yichen Zhao$^{1*}$, and Matthew Hale$^1$
\thanks{*These authors contributed equally.}
\thanks{This work was supported by AFOSR under grant FA9550-19-1-0169, ONR under grant N00014-24-1-2331, and AFRL under grant FA8651-23-F-A006.}}

\maketitle

\begin{abstract}
Asynchronous optimization algorithms often require delay bounds to prove their
convergence, though these bounds can be difficult to obtain in practice.
Therefore, we introduce
a novel distributed generalized momentum algorithm that 
provides fast convergence and 
allows arbitrary finite delays.
It subsumes Nesterov's accelerated gradient
algorithm and the heavy ball algorithm, among others. 
We first develop conditions on the parameters of this algorithm that ensure asymptotic
convergence.
Then we show its convergence rate is linear in a function of the number of computations and communications
that processors perform (in a way that we make precise). 
Simulations compare this algorithm to gradient descent, heavy ball, and Nesterov's accelerated gradient algorithm 
with a text classification problem on the 20 newsgroups dataset. 
Across a range of scenarios with unbounded delays, 
the generalized momentum algorithm converges with at least 36$\%$ fewer iterations than gradient descent, 28$\%$ fewer iterations than the heavy ball algorithm, and 16$\%$ fewer iterations that Nesterov's accelerated gradient algorithm.
\end{abstract}


\section{Introduction}
\label{sec:introduction}
Large-scale convex optimization problems arise in several fields, including robotics~\cite{derenick2007convex,deits2015computing,alonso2015multi}, machine learning~\cite{bottou2018optimization,wang2020survey,bubeck2014theory}, and communications~\cite{mattingley2010real,palomar2010convex,shi2015large}.
Parallel algorithms are commonly used to solve these problems~\cite{upadhyaya2013parallel,salman2023parallel,djafri2022dynamic,verbraeken2020survey}. 
Work in~\cite{bertsekas1989parallel} 
established three categories
for parallelized algorithms: 
(i) ``synchronous'', in which all processors compute and communicate at every iteration, 
(ii) ``partially asynchronous'', in which all processors must compute and communicate at least once
within intervals of a specified length, and
(iii) ``totally asynchronous'', in which there may be arbitrarily long delays between
successive computations and communications performed by each processor. 

Synchronous algorithms are often slowed by the ``straggler penalty,'' which results from processors idling while waiting for the slowest
processor to finish computing and/or communicating~\cite{deshmukh2021collaborative,bertsekas1989parallel}.
Partially asynchronous algorithms can avoid the straggler penalty, but they assume both the
existence and knowledge of a bound on all delays, which
does not always hold. 
Totally asynchronous algorithms are more flexible because they permit arbitrarily long delays in 
communications and computations. However, existing totally
asynchronous algorithms have largely used gradient descent (GD)~\cite{wu2023}
which can converge slowly. 

In this work, we therefore develop a totally asynchronous generalized momentum (GM) algorithm. 
Total asynchrony allows processors to have arbitrarily long delays 
(i) between successive computations, 
(ii) between computing a new iterate and communicating it to other processors,
and 
(iii) between when an iterate is sent
to another processor and when it is received by that processor. 
The GM algorithm
includes the heavy ball (HB) algorithm~\cite{polyak1964some} and Nesterov's accelerated
gradient (NAG) algorithm~\cite{nesterov1983method} as special cases. 

We apply the techniques of~\cite[Chapter 6]{bertsekas1989parallel} to show that 
the GM algorithm converges linearly under total asynchrony.
This approach shows that an update law converges asymptotically under
total asynchrony if it is an $\infty-$norm 
contraction. The GM update law is generally not
an~$\infty$-norm contraction, but 
we show that applying it
twice is. 
We therefore present a two-step GM algorithm and prove
its convergence. It was shown in~\cite{bertsekas1989parallel} that
convergence under total asynchrony requires some form
of diagonal dominance of the Hessian of the objective function.
We likewise enforce a diagonal dominance condition
and our analysis shows how our choice of algorithm
interacts with this property to ensure convergence. 

\subsection{Summary of Main Contributions}
To summarize, our contributions are:
\begin{itemize}
    \item A novel totally asynchronous momentum algorithm guaranteed to converge under arbitrarily long delays in computations and communications (Algorithm \ref{alg:TAGM}).
    \item Bounds on algorithm parameters that ensure convergence is linear in a certain function of the number of 
    computations and communications that processors complete (Theorems~\ref{thrm:sssConv}-\ref{thrm:TAGM}).
    \item Bounds on the minimum number of computations and communications 
    that must be performed to converge within a given distance of an optimum (Theorem \ref{thrm:opcomplex}).
    \item Numerical results that show Algorithm \ref{alg:TAGM} requires at least 36$\%$ fewer iterations than GD, 28$\%$ fewer iterations than HB, and 16$\%$ fewer iterations than NAG on a text classification task (Section \ref{sec:simulations}).
\end{itemize}

\subsection{Related Work}
There is a large body of literature on partially asynchronous optimization algorithms, and a representative sample includes \cite{tian2020asynchronous,zhou2018distributed,assran2020,chang2016,tseng1991rate,tseng1990network,tian2018linear,nedic2001,hsieh2022,sun2017real,magnusson2020,tsianos2011}. 
As noted above, the convergence of these algorithms requires knowledge
of a fixed bound on all delays in all processors' computations and communications. 
Such a bound may not be known or even exist (e.g., if delays grow over time), which 
precludes the use of those algorithms.
We focus on total asynchrony because it does not require the existence of such a delay bound to guarantee the convergence of an algorithm.

Existing totally asynchronous algorithms include 
a projected GD algorithm \cite{bertsekas1989parallel}, the HB algorithm \cite{hustig2023totally}, and 
the NAG algorithm \cite{pond2024}. 
Algorithm~\ref{alg:TAGM} in this paper subsumes all of these algorithms as special cases, and
we empirically show that Algorithm~\ref{alg:TAGM} converges faster
than all of them. 
Additional work in~\cite{mansoori2019fast} presents an asynchronous algorithm
based on Newton's method
in which consecutive computations can have arbitrary delays
between them. However, processors' communications are assumed
to be received immediately in that work, while 
we allow arbitrarily long delays between when a message is sent
and when it is received. 
The algorithm in~\cite{mansoori2019fast} therefore lacks convergence guarantees under
the conditions that we consider and we do not empirically compare to it. 
A similar statement applies to~\cite{srivastava2011}. 

In the centralized setting,
it has been widely shown that momentum algorithms converge faster than gradient 
descent~\cite{haji2021comparison,qian1999momentum,vanscoy2018,drori2014performance,su2016,saab2022adaptive,sun2019}. 
We bring this same acceleration to the totally asynchronous setting for the first time, and empirically demonstrate the computational speed-up obtained when using our proposed algorithm compared to existing algorithms that do not impose delay bounds.

\section{Preliminaries and Problem Statements}
This section defines notation, reviews the generalized momentum algorithm, and states the problems 
we solve.

\subsection{Notation}

The real, non-negative real, and natural numbers are denoted by $\R$, $\R_{\geq0}$, and $\mathbb{N}$ respectively.
We use $\Pi_{\X}:\R^n\rightarrow\X$ to denote the orthogonal projection onto a non-empty, closed, convex set $\X\subset\R^n$, i.e., $\Pi_\X(y) = \argmin_{x\in \X}\Vert x-y\Vert_2$.
For~$x \in \R^n$, the $\infty$-norm is~$\Vert x \Vert_\infty = \max_{i \in \{1, \ldots, n\}}\vert x_i\vert$. 
The Frobenius norm is defined as $\Vert A\Vert_F = \text{trace}(A^TA)^{1/2}$ for a matrix $A\in\R^{m\times n}$.
We use $\vert \cdot \vert$ to denote the cardinality of a set.
The column operator for $a, b \in \R^n$ is $\col{a}{b}= \begin{bmatrix}
a^T & b^T \end{bmatrix}^T\in \mathbb{R}^{2n}$.
For~$f : \R^n \to \R$, we define
$\nabla_if = \frac{\partial f}{\partial x_i}$.

\subsection{The Generalized Momentum Algorithm}
Given an objective function $f:\R^n\rightarrow \R$ and a constraint set $\X\subset \R^n$, we consider the optimization problem 
$\operatorname{minimize}_{x \in \X} \,\, f(x)$.
This problem can be solved by the centralized generalized momentum algorithm~\cite{cyrus2018}: 
at timestep $l\in\mathbb{N}$ the decision variable $x(l)\in \R^n$ is updated according to 
\begin{multline}\label{eq:centralGMupdate}
    x(l+1) = \Pi_\X \Big[ x(l) - \gamma \nabla f\Big(x(l)+ \lambda\big(x(l)-x(l-1)\big) \Big) 
    \\  + \beta\big(x(l)-x(l-1)\big) \Big],
\end{multline}
where $\gamma,\lambda,\beta\in \R_{\geq0}$.
We call it the ``generalized'' momentum method 
because certain choices of
$\gamma$, $\lambda$, and $\beta$
give other momentum algorithms, including Nesterov's accelerated gradient (NAG) algorithm~\cite{nesterov1983method}
and Polyak's heavy ball (HB) algorithm~\cite{polyak1964some}. 
The GM algorithm is thus quite general
and we develop a totally asynchronous version.

\subsection{Problem Statements}
We consider~$n$ processors indexed by the set~$\V = \{1,2,\dots,n\}$ 
and communicating over a graph~$G = (\V, \E)$, where~$\E \subseteq \V \times \V$. 
We use~$\V_i$ to denote processor~$i$'s neighborhood set. 
The decision variable $x$ is partitioned into blocks among the processors, 
with each processor computing new values only of their block of the decision variable.
To ease exposition, we take $x \in \R^n$ to be partitioned into scalar blocks with
$x_i \in \R^{n_i}$ and~$n_i = 1$ for all~$i \in \V$. All analytical
results easily extend to vector blocks with~$n_i > 1$. 
    
We consider objective functions $f:\R^n\rightarrow \R$ of the form
\begin{equation} \label{eq:objfunc}
    f(x) := \sum_{i=1}^n f_i(x_{\V_i}),
\end{equation}
where $f_i:\R^{\vert \V_i\vert+1}\rightarrow \R$ and $x_{\V_i}$ is the vector containing processor $i$'s own decision variable and all decision 
variables of its neighbors, namely processors with indices~$j \in \V_i$.

This paper solves the following three problems. 

\begin{problem}\label{prob:1}
    Construct a totally asynchronous generalized momentum algorithm that solves
    \begin{equation} \label{eq:prob1}
        \underset{x \in \X}{\operatorname{minimize}} \,\, f(x) = \sum_{i=1}^nf_i(x_{\V_i}). 
    \end{equation}    
\end{problem}

\begin{problem}\label{prob:2}
    Show that the totally asynchronous GM algorithm in Problem \ref{prob:1} converges linearly to a minimizer.
\end{problem}

\begin{problem}\label{prob:3}
    Given $\epsilon > 0$, find the number of computations and communications each processor must execute 
    such that the iterates of the totally asynchronous GM algorithm are within a distance $\epsilon$ of the solution to~\eqref{eq:prob1}. 
\end{problem}

with \section{Totally Asynchronous Distributed Generalized Momentum Method}

In this section we solve Problem \ref{prob:1} by formulating a totally asynchronous GM algorithm and lay the groundwork for solving Problem \ref{prob:2}.
We apply the framework presented in \cite{bertsekas1989parallel} to establish the convergence of a totally asynchronous algorithm based on the properties of the synchronous version of the algorithm.
We will first show that two variations of the synchronous GM algorithm satisfy certain technical conditions in order to present the totally asynchronous version of the algorithm.
Specifically, we will show that two applications of the GM method result in an $\infty$-norm contraction map, which allows us to guarantee asymptotic convergence of the totally asynchronous version.
We emphasize that all discussions in this section regarding the two synchronous algorithms
are steps toward proving the convergence of the totally asynchronous GM algorithm.

We make the following three assumptions about the optimization problem in~\eqref{eq:prob1}.

\begin{assumption}\label{ass:X}
    The constraint set $\X\subset \R^n$ is nonempty, convex, and compact. The set can be decomposed as $\X=\X_1 \times \X_2 \times \dots \times \X_n$, where $\X_i\subset \R^{n_i}$ for all $i\in \V$.
\end{assumption}

This assumption on $\X$ is necessary for parallelizing the projection operator in the centralized GM update law in~\eqref{eq:centralGMupdate}, which
we will use to 
guarantee constant satisfaction of the set constraint in \eqref{eq:prob1}, even under total asynchrony.

\begin{assumption}\label{ass:f}
    The objective function $f$ is twice continuously differentiable.
\end{assumption}

Assumption \ref{ass:f} guarantees the existence and continuity of the gradient and Hessian of $f$, which we 
use in the convergence analysis of the totally asynchronous GM algorithm.

\begin{assumption} \label{ass:H}
    The Hessian matrix~$H(x):=\nabla^2f(x)\in \R^{n\times n}$ is $\mu-$diagonally dominant on $\X\subset \R^n$ for some $\mu>0$. 
    That is, for each $i\in \V$ we have $H_{ii}(x)\geq \mu + \sum_{j=1,j\not=i}^n\vert H_{ij}(x) \vert$ for all $x\in \X$.
\end{assumption}

It was shown in~\cite[Section 6.3.2]{bertsekas1989parallel} that 
    some form of Hessian diagonal dominance is 
    generally needed to establish convergence of an algorithm under total asynchrony. 
    In fact, \cite[p. 438]{bertsekas1989parallel} explicitly showed that
    a lack of diagonal dominance can lead to divergence. 
    Therefore, we enforce Assumption \ref{ass:H}. 
    Intuitively, it asserts that the gradient~$\nabla_i f(x)$
    depends more on~$x_i$ than on all other entries of~$x$. 
    Settings where diagonal dominance appears include text classification \cite{greene2006practical}, bioinformatics \cite{scholkopf2002kernel}, sums-of-squares optimization \cite{ahmadi2015sum}, iterative algorithms for solving large linear systems \cite{cai2019fast,gremban1996combinatorial}, 
    and quadratic programs \cite{ubl2021totally}.
    Generally, it was noted by \cite{greene2006practical,cancedda2003word,scholkopf2002kernel} that diagonal dominance is a common structural property of kernel methods applied to sparse, high-dimensional data \cite{greene2006practical,cancedda2003word,scholkopf2002kernel}.
    For problems where diagonal dominance is not inherent, it can often be imposed through regularization \cite{hendrickson2022totally,koshal2011multiuser}.

Assumption \ref{ass:H} implies that $f$ is $\mu-$strongly convex on $\X$
and thus it has a unique minimizer over $\X$, which is denoted by
\begin{equation} \label{eq:xstardef}
x^\star=\col{x_1^\star,\dots}{x_n^\star}.
\end{equation}

The GM update law in~\eqref{eq:centralGMupdate}
depends on the iterate $x(l)$ and the prior iterate $x(l-1)$. 
Throughout this paper, we use $y(l) = x(l-1)$.
For each $i\in \V$, processor~$i$ will store a local copy of both~$x(l)$ and~$y(l)$ onboard, 
which we denote $z^i(l):=(x^i(l),y^i(l))\in \Z := \X\times\X$, where
the superscript~$i$ indicates a vector that is stored onboard processor~$i$.
\begin{definition}\label{def:zstar}    
    We define $z^\star=(x^\star,x^\star)\in \Z$ and $z^\star_i=(x^\star_i,x^\star_i)\in \Z_i := \X_i \times \X_i$, where
    $x^\star$ is from~\eqref{eq:xstardef}. 
\end{definition}

\begin{table}[t]
\centering
\begin{tabular}{lll}
        \toprule
        Symbol & Dimension &
        Description \\
        \midrule
        \multirow{2}{*}{$z^i=(x^i,y^i)$ }& \multirow{2}{*}{$\R^n\times \R^n$} & Processor $i$'s local copy of the decision\\ &&  vector \\
        \hline
        \multirow{2}{*}{$x^i(l)$ } & \multirow{2}{*}{$\R^n$} & Processor $i$'s local copy of $x$ at \\ && time step $l$ \\ \hline
        \multirow{2}{*}{$y^i(l)$}& \multirow{2}{*}{$\R^n$}& Processor $i$'s local copy of $y$ at \\ & & time step $l-1$ \\\hline
        \addlinespace
        $x^i_i,y^i_i$ & $\R^{n_i}$ &  Processor $i$'s local copy of its own\\ && decision variables \\
        \hline
        \addlinespace $x^i_j,y^i_j$ & $\R^{n_j}$  & Processor $i$'s local copy of processor\\ &&  $j$'s decision variables \\ \addlinespace \hline \addlinespace
        $\X_i$ & $\R^{n_i}$ & Constraint set for $x^i_i$ 
        \\ \addlinespace \hline \addlinespace
        $\Z_i$ & $\R^{n_i\times n_i}$ & Constraint set for $z_i^i$ \\ \addlinespace
        \bottomrule
\end{tabular}
\caption{
Superscripts indicate which processor stores the decision vector, while subscripts index the entries of that vector.
}
\label{tab:notation}
\end{table}


Over time, processor~$i$ will compute new values of its
block of decision variables, 
which in this notation is~$z^i_i(l)=(x^i_i(l),y^i_i(l))\in \Z_i$, where 
the subscripts index entries of each vector. 
The terms~$z^i_i$, $x^i_i$, and~$y^i_i$ indicate that processor $i$ is responsible for updating these variables and that they are stored onboard processor $i$, 
whereas $z^i_j$, $x^i_j$, and $y^i_j$ with~$j \neq i$ are stored onboard processor $i$ but processor $j$ is responsible for updating these variables and communicating
new values of them to processor~$i$. 
Table \ref{tab:notation} summarizes this notation.
At time~$l$ processor~$i$ stores onboard the quantities 
\begin{multline}
    z^i(l) = \big( \col{x^i_1(l),\dots,x^i_i(l),\dots}{x^i_n(l)},
    \\ \col{y^i_1(l),\dots,y^i_i(l),\dots}{y^i_n(l)}\big).
\end{multline}
The objective~$f$ in~\eqref{eq:objfunc} shows that~$f_i$
depends only on the decision variables of processor~$i$
and its neighbors. 
Only those decision variables appear
when processor~$i$ computes
the gradient~$\nabla_i f$ and in particular no decision variable~$x^i_m$
or~$y^i_m$ appears if~$m \not\in \V_i$. 
Processor $i$ can therefore set its copy of~$x^i_m$ and $y^i_m$ to arbitrary values for  
processors with indices $m \not\in \V_i$. 
Those decision variables
do not affect its computations and 
processors~$i$ and~$m$ never communicate. 

\subsection{The Single-Step Synchronous Method}
In this sub-section, we define and analyze the first variation of the synchronous GM algorithm, which we refer to as the ``single-step'' 
synchronous method, which we later use to create the totally asynchronous GM algorithm.

The single-step synchronous GM (sGM) algorithm has simultaneous computations and communications among processors. 
Processor~$i$ updates~$x^i_i$ and $y^i_i$ using the maps $\tilde{u}^i_x,\tilde{u}^i_y:\Z\rightarrow \X_i$ defined as
\begin{align}
    x^i_i(l+1) &= \tilde{u}^i_x\big(x^i(l),y^i(l)\big)  \label{eq:sssuxi}
    \\ &= \Pi_{\X_i}\Big[x^i_i(l) - \gamma\nabla_i f\Big(x^i(l)+\lambda\big(x^i(l)-y^i(l)\big)\Big) 
    \\ &\qquad\qquad\qquad\qquad\qquad + \beta(x^i_i(l)-y^i_i(l))\Big]
    \\y^i_i(l+1) &= \tilde{u}^i_y\big(x^i(l),y^i(l)\big) \label{eq:sssuyi}
    \\ &= x^i_i(l)
\end{align}
for all $l\in\mathbb{N}$ and $i\in \V$.
Let $\tilde{u}^i:\Z\rightarrow \Z_i$ 
denote the combined update map
\begin{align} \label{eq:sssui}
    \big(x^i_i(l+1)&,y^i_i(l+1)\big) = \tilde{u}^i\big(x^i(l),y^i(l)\big)
    \\ &= \Big(\tilde{u}^i_x\big(x^i(l),y^i(l)\big), \tilde{u}^i_y\big(x^i(l),y^i(l)\big)\Big),
\end{align}
which we will also write as~$z^i_i(l+1) = \tilde{u}^i(z^i(l))$. 
This update law performs one iteration of GM and stores the result in~$x^i_i(l+1)$ in \eqref{eq:sssuxi}; 
the variable~$y^i_i(l+1)$ stores the previous value of~$x^i_i$, which is~$x^i_i(l)$. 

We next show the sGM update law in~\eqref{eq:sssui} has a unique fixed point, which is the solution to~\eqref{eq:prob1}.

\begin{lemma}\label{lem:sssFP}
    Consider Problem \ref{prob:1} and let Assumptions~\ref{ass:X}-\ref{ass:H} hold. 
    Then the point $z^\star$ from Definition~\ref{def:zstar}
    is a fixed point of the single-step synchronous GM update law \eqref{eq:sssui}, 
    in the sense that $z^\star_i = \tilde{u}^i(z^\star)$ for all $i\in \V$.
\end{lemma}
\begin{proof}
    See Appendix~\ref{app:lem1proof}. 
\end{proof}

For synchronous algorithms, every processor updates its decision variables and communicates at each timestep $l\in\mathbb{N}$.
Processor $i\in\V$ updates according to \eqref{eq:sssui} and communicates $z^i_i(l+1)=(x^i_i(l+1),y^i_i(l+1))$ to all 
processors~$j \in \V_i$.
Also at timestep $l\in\mathbb{N}$, all processors~$j \in \V_i$ use this communication from processor $i$ 
to overwrite the prior value of processor~$i$'s decision variables that processor~$j$ has onboard.
Mathematically, processor~$j$ sets~$z^j_i(l) = z^i_i(l)$ for each~$j \in \V_i$.

We define the true state of the network 
 at time $l\in\mathbb{N}$ to be the vector that collects the current value of each decision variable from each processor.
 Formally,~$z^\tr(l)= \big(x^\tr(l),y^\tr(l)\big)$, where~$x^\tr(l)=\col{x^1_1(l),\dots}{x^n_n(l)}$ and~$y^\tr(l) =\col{y^1_1(l),\dots}{y^n_n(l)}$.
In the next theorem, we show that the sGM algorithm is an~$\infty$-norm contraction mapping when it is applied
twice to the true state of the network.
In it, we refer to the map 
\begin{multline}\label{eq:utrue}
    \hat{u}^i_\tr\big(z^\tr(l)\big) = x^i_i(l)-\gamma\nabla_i f\big(x^\tr(l)\\+\lambda\big(x^\tr(l)-y^\tr(l)\big)\big) + \beta\big(x^i_i(l)-y^i_i(l)\big),
\end{multline}
which models the evolution of the~$i^{th}$ entry of the true state of the network at each iteration $l\in\mathbb{N}$.
The update law in~\eqref{eq:utrue} 
only holds in the synchronous setting
because only then does each processor have all of its neighbors' most recent iterates. 
The analysis of the synchronous algorithm will enable us to derive
conditions for convergence of its totally asynchronous counterpart
in the next sub-section.
Given~$\mu > 0$ from Assumption~\ref{ass:H}, 
the forthcoming theorem uses the sets
    \begin{multline}\label{eq:C1}
        C_1=\bigg\{(\gamma,\lambda,\beta)\in\R\times\R\times \R : 
        0 \leq \beta\leq \lambda< \frac{\gamma\mu}{2(1-\gamma\mu)}, 
        \\0<\gamma < \frac{\beta}{\lambda \max\limits_{i\in\V}\max\limits_{\eta\in\X}\vert H_{ii}(\eta)\vert} \bigg\}
        \end{multline}
    and
    \begin{multline}\label{eq:C2}
        C_2 = \bigg\{ (\gamma,\lambda,\beta)\in\R\times\R\times\R : 0\leq \lambda \leq \beta< \frac{1}{2}\gamma\mu(1+2\lambda), 
        \\ 0 <\gamma< \frac{1}{\max\limits_{i\in\V}\max\limits_{\eta\in\X}\vert H_{ii}(\eta)\vert}\Bigg\}.
    \end{multline}
The following theorem also 
uses the constants
\begin{align}\label{eq:a1}
    \alpha_1&= \big(1+\beta-\gamma\mu(1+\lambda)\big)^2 
    \\ &\quad\quad\quad + (\beta-\gamma\lambda \mu)\big(2+\beta-\gamma\mu(1+\lambda)\big)
\end{align}
and
\begin{align}\label{eq:a2}
        \alpha_2 = 1-\gamma\mu+2(\beta-\lambda\gamma\mu).
\end{align}

The contraction properties of the sGM update law are as follows. 

\begin{theorem}\label{thrm:sssConv}
    Consider the problem in~\eqref{eq:prob1}, let Assumptions~\ref{ass:X}-\ref{ass:H} hold, and
    let~$z(0) \in \Z$ be given. 
    Consider the sGM update law in~\eqref{eq:sssui}. 
    If the parameters $(\gamma,\lambda,\beta)\in C_1\cup C_2$,
    then the sGM update law satisfies~$\Vert z(l+1)-z^\star\Vert_\infty\leq \alpha\Vert z(l-1)-z^\star\Vert_\infty$ for all $l\in\mathbb{N}$, where $\alpha=\max\{\alpha_1,\alpha_2\}\in(0,1)$ with $\alpha_1$ defined in \eqref{eq:a1} and $\alpha_2$ defined in \eqref{eq:a2}.
\end{theorem}
\begin{proof}
See Appendix~\ref{app:main}. 
\end{proof}

Theorem \ref{thrm:sssConv} proves that the sGM algorithm is contractive with respect to the $\infty$-norm over two timesteps, 
i.e., it is contractive from timestep~$l-1$ to timestep~$l+1$.
We use this result in the next section
to create an alternate synchronous GM algorithm that is contractive over one timestep.

\subsection{The Double-Step Synchronous Method}
We now develop a synchronous GM algorithm that is contractive over one timestep in order to 
apply the method in~\cite{bertsekas1989parallel} to prove convergence under total asynchrony. 
The synchronous algorithm that we develop in this sub-section is termed the ``double-step synchronous GM'' (dGM) algorithm
because each processor computes two applications of the GM update law~\eqref{eq:centralGMupdate} at each timestep.
After we define the dGM algorithm, we show that it satisfies a three-part lemma from \cite{bertsekas1989parallel} 
that guarantees asymptotic convergence of its totally asynchronous form.

From now on, we use the variable $k\in\mathbb N$ to represent time in order to
distinguish the sGM algorithm from the dGM algorithm.
Incrementing~$k$ to~$k+1$ is equivalent to incrementing~$l$ to~$l+2$. 
Under the dGM algorithm, each processor updates its decision variable $z^i_i=(x^i_i,y^i_i)$ at every $k\in\mathbb{N}$ with the maps $u^i_x,u^i_y:\Z\rightarrow \X_i$, defined as
\begin{align}
    x^i_i&(k+1) = u^i_x\big(x^i(k),y^i(k)\big) \label{eq:uix}
    \\ &= \Pi_{\X_i}\Big[ y^i_i(k+1)+ \beta\big(y^i_i(k+1)-x^i_i(k)\big)
    \\&\quad-\gamma\nabla_i f\Big( y^i(k+1) 
    + \lambda\big(y^i(k+1)-x^i(k)\big)\Big) \Big]
   \\ y^i_i&(k+1) = u^i_y\big(x^i(k),y^i(k)\big)\label{eq:uiy}
   \\ &= \Pi_{\X_i}\Big[ x^i_i(k) + \beta\big(x^i_i(k)-y^i_i(k)\big)
   \\&\quad - \gamma\nabla_i f\Big( x^i(k)+ \lambda\big(x^i(k)-y^i(k)\big)\Big)\Big].
\end{align}
In \eqref{eq:uix}, we define $y^i(k+1)$ as $y^i(k+1)=\col{y^i_1(k),\dots,y^i_i(k+1),\dots}{y^i_n(k)}$, where $y^i_i(k+1) \in \mathbb{R}$ is the newly updated local variable in \eqref{eq:uiy} and all other entries satisfy~$y^i_j(k+1) = y^i_j(k)$ 
for $j \in \V\backslash\{i\}$.
The term~$y^i_i(k+1)$ in~\eqref{eq:uix} can be expanded so that~\eqref{eq:uix}
only depends on~$z^i(k)$ and not any iterates from time~$k+1$. 
Then~\eqref{eq:uix} and~\eqref{eq:uiy} can be computed in a single timestep.

%

Let $u^i:\Z\rightarrow\Z_i$ be the combined update map for the dGM algorithm so 
that $z^i_i(k+1) = u^i(z^i(k))= ( u^i_x(z^i(k)),u^i_y(z^i(k)))$ holds for all $k\in\mathbb N$.

\begin{lemma}\label{lem:dGMfp}
    Consider the problem in~\eqref{eq:prob1}
    and let Assumptions~\ref{ass:X}-\ref{ass:H} hold. 
    Then~$z^\star$ from Definition~\ref{def:zstar}
    is a fixed point of the double-step synchronous GM update law in \eqref{eq:uix}-\eqref{eq:uiy}, in the sense that $z^\star_i=u^i(z^\star)$ for all $i\in\V$.
\end{lemma}
\begin{proof}
    Omitted due to similarity to Lemma \ref{lem:sssFP}.
\end{proof}

Let the map $h:\Z\rightarrow \Z$ be defined as
\begin{multline}\label{eq:h}
    h(z) = \text{col}\Big(\big(u^1_x(z),u^2_x(z) \dots, u^n_x(z)\big)^T,\\ \big(u^1_y(z),u^2_y(z),\dots, u^n_y(z)\big)^T\Big),
\end{multline}
which models all processors' computations at one timestep~$k\in\mathbb N$ under the dGM algorithm.
Lemma~\ref{lem:dGMfp} implies that $z^\star$ is the fixed point of $h$, i.e., $h(z^{\star}) = z^{\star}$. 

The mapping~$h$ is the synchronous double-step GM update law and we can analyze
it to ensure convergence under total asynchrony. 
Our main tool in doing so is a three-part lemma from \cite{bertsekas1989parallel}
that we state here in a slightly stronger form than the original. 

\begin{lemma}[\cite{bertsekas1989parallel}] \label{lem:threepart}
    Consider the problem in~\eqref{eq:prob1} and let Assumptions~\ref{ass:X}-\ref{ass:H} hold.
    Suppose that an initial solution estimate $z(0)\in\Z(0)$ is given.
    An update law~$\chi : \Z \to \Z$ converges under total asynchrony
    if there exist sets~$\{\Z(k)\}_{k \in \mathbb{N}}$ that satisfy: 
    \begin{enumerate}
        \item 
        The Lyapunov-like containment (LLC) condition: the containment relationship
        \begin{equation}
            \Z \supset \Z(0)\supset \cdots \supset\Z(k) \supset\Z(k+1)\supset \cdots
        \end{equation}
        holds for all $k\in\mathbb N$
        \item The synchronous convergence condition (SCC): 
        (i) $z\in\Z(k)$ implies $\chi(z)\in \Z(k+1)$ for all $k\in\mathbb N$ and (ii) a sequence $\{z_k\}_{k\in\mathbb N}$ with $z_k\in\Z(k)$ for all $k\in\mathbb N$ satisfies $\lim\limits_{k\rightarrow \infty}z_k=z^\star$, where $z^\star$ is the fixed point of $\chi$.
        \item The box containment condition (BCC): for all $i\in\V$ there exists~$\Z_i(k)\subset \Z_i$ such that $\Z(k)=\Z_1(k)\times \cdots \times \Z_n(k)$.
    \end{enumerate}
\end{lemma}

We next use the result of Theorem~\ref{thrm:sssConv} to 
show that Lemma~\ref{lem:threepart} is satisfied and hence
establish that dGM asymptotically converges under
total asynchrony. 

\begin{theorem} \label{thrm:asymptotic}
Consider the problem in~\eqref{eq:prob1} and let Assumptions~\ref{ass:X}-\ref{ass:H} hold. 
Let~$h$ be the dGM update law from~\eqref{eq:h}. 
Let an initial estimate~$z(0) \in \Z(0)$ be given and
define the sets~$\{\Z(k)\}_{k \in \mathbb{N}}$ as
    \begin{equation}\label{eq:Zk}
        \Z(k) = \{v\in \Z : \Vert v-z^\star\Vert_\infty \leq \alpha^k\Vert z(0)-z^\star\Vert_\infty\},
    \end{equation}
    where~$\alpha$ is from Theorem~\ref{thrm:sssConv}. 
    Then these sets satisfy all conditions in Lemma~\ref{lem:threepart} and hence~$h$
    converges asymptotically under total asynchrony. 
\end{theorem}

\begin{proof}
    See Appendix~\ref{app:zsets}. 
\end{proof}

\subsection{The Totally Asynchronous GM Algorithm}
We now use the synchronous sGM and dGM algorithms to develop the totally asynchronous GM algorithm. 
All processors use the same update law $u^i$ as in the dGM algorithm defined in \eqref{eq:uix}-\eqref{eq:uiy}.
Total asynchrony means the processors do not have to compute and communicate at every timestep $k\in\mathbb N$.
Let $K^i\subseteq \mathbb N$ be the set of times at which processor $i \in \V$ performs a computation.
If $k\in K^i$, then $(x^i_i(k+1),y^i_i(k+1))\gets u^i(x^i(k),y^i(k))$.
If~$k \not\in K^i$, then~$(x^i_i(k+1),y^i_i(k+1)) \gets (x^i_i(k),y^i_i(k))$, i.e., 
processor~$i$'s decision variables are held constant. 
The sets~$\{K^i\}_{i \in \V}$ are solely used for analysis and processors do not need to know them to execute
the forthcoming algorithm. 

Let the set $R^i_j\subseteq\mathbb N$ be the set of times at which processor $i$ receives a communication from processor $j$. 
Because we allow for arbitrarily long communication delays, processor~$i$ may receive such communications many timesteps after they are sent. 
If $k\in R^i_j$, then $(x^i_j(k),y^i_j(k))\gets ( x^j_j(\tau^i_j(k)),y^j_j(\tau^i_j(k)))\in\Z_j$, where 
the map~$\tau^i_j : \mathbb{N} \to K^j$ outputs the 
time at which processor $j$ originally computed the decision variable that processor~$i$ has onboard at time~$k$.
That is, $x^i_j(k) = x^j_j(\tau^i_j(k))$ for all~$k \in \mathbb{N}$. 
As with~$K^i$, the sets~$\{R^i_j\}_{i, j \in \V}$ are only used for ease of exposition and are not known by the processors. 
If~$m \not\in \V_i$, then communications never occur between processors~$i$ and~$m$
and hence~$R^i_m=\emptyset$; 
the value of~$(x^i_m(k),y^i_m(k))$ can be arbitrary at all times~$k \in \mathbb{N}$
because it does not affect processor~$i$'s computations.

The totally asynchronous GM algorithm is formalized using these rules in Algorithm \ref{alg:TAGM}, which solves
Problem~\ref{prob:1}.
We note that a property of Algorithm \ref{alg:TAGM} is that processors receive communications in the order in which they were sent.

\begin{algorithm}
    \caption{Totally Asynchronous GM Algorithm}
    \label{alg:TAGM}
	\SetAlgoLined
	\KwIn{For $i\in\mathcal{V}$ select an initial state $z^i(0)\in\Z.$}
	\For{$k\in\mathbb{N}$}
	{
		\For{$i\in\mathcal{V}$}
		{
			\If{$k\in K^i$}
            { ~\phantom{} $\big(x^i_i(k+1),y^i_i(k+1)\big)\gets u^i\big(x^i(k),y^i(k)\big)$

                \If{$j\in\mathcal{V}^i$}
                { ~\phantom{} Send $\big(x^i_i(k+1),y^i_i(k+1)\big)$ to processor $j$
                }   
			}			
            \For{$j\in\mathcal{V}^i$}
            {
                \If{$k\in R^i_j$}
                { ~\phantom{} $\big(x^i_j(k),y^i_j(k)\big)\gets \Big(x^j_j\big(\tau^i_j(k)\big),y^j_j\big(\tau^i_j(k)\big)\Big)$
                }
            }
            \If{$m\not\in\mathcal{V}^i$}
            { ~\phantom{} $\big(x^i_m(k),y^i_m(k)\big)\gets \big(x^i_m(k-1),y^i_m(k-1)\big)$
            }
		}	
	}
\end{algorithm}

We reiterate that (i) the set~$K^i$ does not need to be known by any processor,
(ii) the set~$R^i_j$ does not need to be known by any processor,
and (iii) in line~$6$ processor~$j$ can receive communications from processor~$i$
after arbitrarily long delays, i.e., communications need not
be received immediately. 

\section{Convergence Analysis}
This section solves Problem \ref{prob:2} and establishes that the totally asynchronous GM algorithm in Algorithm \ref{alg:TAGM} converges linearly.
First, we state a standard 
assumption. 

\begin{assumption}[\cite{bertsekas1989parallel}]\label{ass:TA}
    The sets $K^i$ and $R^i_j$ are infinite for all processors~$i \in \V$ and~$j \in \V_i$.
    Additionally, if $\{k_s\}_{s\in \mathbb{N}}$ is an increasing sequence of times in $K^i$, then
    $\lim\limits_{s\rightarrow\infty}\tau^i_j(k_s)=\infty$ for all $i\in\V$ and $j\in\V_i$.
\end{assumption}

Assumption \ref{ass:TA} guarantees that no processor will ever permanently 
stop computing or communicating, though it allows arbitrarily long delays between these events.
We emphasize that although each individual delay has finite length, the delays themselves can be arbitrarily long, and they can even grow without bound over time. 
For example, under the totally asynchronous model, a processor may wait~$2^k$ timesteps
before sending its~$k^{th}$ communication to a neighbor. 
Over time, this delay grows arbitrarily large and eventually exceeds any fixed bound. 
Such a pattern of delays is permitted by our algorithm.

To analyze the convergence properties of Algorithm \ref{alg:TAGM}, we define an ``operation cycle'' 
using the map~$\text{ops}:\mathbb{N}\rightarrow \mathbb{N}$, 
which counts communications and computations in the totally asynchronous GM algorithm in the following
way.
When $k=0$, we have $\ops k=0$.
Let $k'\in\mathbb N$ be the first point in time when, for all~$i \in \V$, 
(i)~processor~$i$ has updated its decision variables, 
(ii)~processor~$i$ has sent its updated decision variables to all $j\in\V_i$, 
and 
(iii)~for all~$j \in \V_i$, processor~$j$
has received these communications. 
At this point $k'\in\mathbb N$, we set $\ops {k'}=1$ to reflect the completion of a full cycle of Algorithm \ref{alg:TAGM}.
For all $k\geq k'$, we have $\ops{k'}=1$ until the time $k'' \in \mathbb{N}$ 
at which all steps (i)-(iii) have been completed again, at which point $\ops{k''} = 2$.
Assumption~\ref{ass:TA} implies that
$\lim\limits_{k\rightarrow\infty}\ops{k}=\infty$.

Processors may compute and communicate more than once per operation cycle.
However, the value of~$\text{ops}(k)$
is not incremented until \textit{every} processor has computed, sent, and received data. 
At $\ops{k}=0$, each processors' local decision variable obeys $(x^i(k),y^i(k))\in\Z(0)$.
At time $k'$ with $\ops{k'}=1$, we have $(x^i(k'),y^i(k'))\in \Z(1)$ for all $i\in\V$.
We next establish an invariance property of the sets $\Z(k)$ that is then used to derive
a convergence rate. 

\begin{lemma} \label{lem:invariant}
    Consider using Algorithm~\ref{alg:TAGM} on the problem in~\eqref{eq:prob1}
    and let Assumptions~\ref{ass:X}-\ref{ass:TA} hold. 
    Let initial estimates~$z^i(0) \in \Z(0)$ for all $i\in\V$ be given.
    Then each set $\Z(k)$ defined in \eqref{eq:Zk} is forward invariant under Algorithm~\ref{alg:TAGM}, i.e., 
    once $z^i(k_1)\in\Z(k)$ for all $i\in\V$ and some $k_1\in\mathbb N$, then $z^i(k_2)\in\Z(k)$ for all 
    $i \in \V$ and all~$k_2\geq k_1$.
\end{lemma}
\begin{proof}
    See Appendix~\ref{app:invariant}. 
\end{proof}


In the following theorem, we solve Problem \ref{prob:2} and prove that the totally asynchronous GM algorithm converges linearly in the value of $\ops{k}$.

\begin{theorem}\label{thrm:TAGM}
    Consider using Algorithm~\ref{alg:TAGM} on the problem in~\eqref{eq:prob1} and 
    suppose that Assumptions~\ref{ass:X}-\ref{ass:TA} are satisfied.
    Let initial estimates $z^i(0)\in\Z(0)$ for all $i\in\V$ be given.
    For each $(\gamma,\lambda,\beta)\in C_1\cup C_2$, 
    the iterates of Algorithm~\ref{alg:TAGM} satisfy
    \begin{equation}
        \max_{i\in\V}\Vert z^i(k)-z^\star\Vert_\infty \leq \alpha^{\ops{k}}\max_{i\in\V}\Vert z^i(0)-z^\star\Vert_\infty
    \end{equation}
    for all $k\in\mathbb N$, where $\alpha\in(0,1)$ is from Theorem \ref{thrm:sssConv}.
\end{theorem}
\begin{proof}
    See Appendix~\ref{app:linear}. 
\end{proof}

The following corollary establishes the relationship between the totally asynchronous GM algorithm and other totally asynchronous gradient algorithms in the literature.

\begin{corollary}\label{corr:comparison}
    Consider using Algorithm~\ref{alg:TAGM} on the problem in~\eqref{eq:prob1} and 
    suppose that Assumptions~\ref{ass:X}-\ref{ass:TA} are satisfied.
    Then the following three conditions hold:
    \begin{enumerate}
        \item If $(\gamma,\lambda,\beta)\in C_2$ and $\beta=\lambda=0$, then the totally asynchronous GM algorithm is equivalent to the totally asynchronous projected GD algorithm in \cite{bertsekas1989parallel}.
        \item If $(\gamma,\lambda,\beta)\in C_2$ and $\lambda = 0$, then the totally asynchronous GM algorithm is equivalent to the totally asynchronous HB algorithm in \cite{hustig2023totally}.
        \item If $(\gamma,\lambda,\beta)\in C_1$ and $\beta=\lambda$, then the totally asynchronous GM algorithm is equivalent to the totally asynchronous NAG algorithm in \cite{pond2024}.
    \end{enumerate}
\end{corollary}

\section{Operation Complexity}
Now we solve Problem~\ref{prob:3} by quantifying the relationship between the 
topology of the processors' network $\G=(\V,\E)$ and the ``operation complexity'' of Algorithm~\ref{alg:TAGM},
which we define to be the number of computations and communications that processors must
execute to compute a solution within a desired accuracy bound. 

Theorem~\ref{thrm:TAGM} established that 
one operation cycle is completed once every processor has 
computed and received information from each of its neighbors. 
Restated in terms of the graph $\G$, one operation cycle consists of at least $\vert \V\vert$ computations 
(one per processor)
and $2\vert \E\vert$ communication events, 
where the factor of two arises because communication must occur both ways
between each pair of neighbors. 

Let $D = \max_{v_1,v_2\in\Z}\Vert v_1-v_2\Vert_\infty$ be the~$\infty$-norm diameter
of the set~$\Z$. 
The following theorem solves Problem \ref{prob:3} by quantifying the minimum number of computations and communications that must be executed by each processor in order for all processors' iterates to be within distance~$\epsilon > 0$ of the minimizer.

\begin{theorem}\label{thrm:opcomplex}
    Consider using Algorithm~\ref{alg:TAGM} on the problem in~\eqref{eq:prob1} and 
    suppose that Assumptions~\ref{ass:X}-\ref{ass:TA} are satisfied.
    Let~$z^i(0) \in \Z(0)$ be given for all~$i \in \V$ and let
    algorithm parameters $(\gamma,\lambda,\beta)\in C_1 \cup C_2$ be given. 
    For all processors' iterates to be within distance $\epsilon > 0$ of the minimizer, i.e., for~$\max\limits_{i\in\V}\Vert z^i(k)-z^\star\Vert_\infty\leq \epsilon$, 
    it is sufficient for processor~$i$ to perform~$\rho$ 
    computations and send~$\rho\vert \V_i\vert$ communications
    for each~$i \in \V$, where $\rho :=\frac{\log (D/\epsilon)}{\log(1/\alpha)}$. 
\end{theorem}
\begin{proof}
    See Appendix~\ref{app:complex}. 
\end{proof}

\section{Simulations} \label{sec:simulations}

    To compare Algorithm \ref{alg:TAGM} to other totally asynchronous algorithms, we solve a distributed learning problem for text classification of the 20newsgroups dataset \cite{newsgroups}, which
    consists of approximately 18,000 documents. Each
    document is about one of $K=20$ different subjects 
    (e.g., space, motorcycles, religion, etc.).
    We use term frequency-inverse document frequency (TF-IDF) to numerically represent the importance of a word to a document.
    Common words that do not carry much meaning, such as ``and'' or ``the'', are not included in the TF-IDF conversion. 
    The dataset is partitioned so that 90\% of the data points are used for training and 10\% are used for testing. 
    The feature dimension of each sample is set to $D=12,000$ so that each feature corresponds to one of the 12,000 most frequently appearing words within the collection of documents.

    We solve the optimization problem
    \begin{equation}\label{eq:newsgroupF}
        \min_{W\in\R^{K\times D}}\; f(W)=\frac{1}{2M}\Vert \Phi W^T- Y\Vert_{F}^2 + \frac{\theta}{2}\Vert W\Vert_F^2,
    \end{equation}
    where $M=16,950$ is the total number of data points
    used for training, 
    $W\in\R^{K\times D}$ is the weight matrix, the rows of $\Phi\in\R^{M\times D}$ are the feature vectors of the form $\phi_i\in\R^D$, the rows of $Y\in\R^{M\times K}$ are one-hot vectors of the form $y_i\in\R^K$ that encode the class label of $\phi_i\in\R^D$, 
    and $\theta>0$ is the $\ell_2$-regularization parameter. 
    The gradient of this function is
    \begin{equation}
        \nabla f(W)= \frac{1}{M}(\Phi W^T - Y)^T\Phi + \theta W,
    \end{equation}
    which we vectorize with the operator $\text{vec}:\R^{K\times D}\rightarrow \R^{KD}$ for numerical computations to form 
    \begin{equation}
        \text{vec}({\nabla f (W)}) = \frac{1}{M}\big((\Phi^T\Phi \otimes I_K)\text{vec}(W) -\text{vec}(Y^T\Phi)\big)+ \theta \text{vec}(W),
    \end{equation}
    where $\text{vec}(W)\in\R^{KD}$.
    We compute the Hessian $\nabla^2(\text{vec}(f(W)))\in\R^{KD\times KD}$ as
    \begin{equation}
        \nabla^2(\text{vec}(f(W))) = \frac{1}{M}\Phi^T\Phi\otimes I_K + \theta I_{KD}.
    \end{equation}
    A straightforward calculation shows that $\nabla^2(\text{vec}(f(W)))$ satisfies Assumption \ref{ass:H} for all $W$ if $\theta\geq 0.16$.
    We set $\theta=0.66$, which makes
    the Hessian of $f$ $\mu$-diagonally dominant with $\mu=0.5$.

    The problem is distributed among $N=50$ processors by assigning $d=D/N$ columns of $W$ to each 
    processor as decision variables to update.
    Processor $n\in[N]$ updates and communicates values of~$W^n_n\in\R^{K\times d}$ to all other processors in the network.
    The underlying training data is distributed among the network by giving $m=M/N$ rows of $\Phi$ (which are training data points) to each processor, along with their labels. 
    Let $\Phi^n\in\R^{m\times D}$ and $Y^n\in\R^{m\times K}$ be the portions of the training data set that processor $n\in[N]$ is given.
    Then the local gradient for processor $n\in[N]$ is
    \begin{equation}
       \nabla_n f(W^n)=\frac{1}{M}(\Phi^n(W^{n})^T-Y^n)^T\Phi^n+\frac{\theta}{N}W^n.
        \end{equation}
    Processor $n\in[N]$ needs access to the full decision variable $W\in\R^{K\times D}$ in order to compute $\Phi^n(W^{n})^T$.
    Therefore, the communication network is the complete graph on $N$ nodes.

    A total of 8 simulations were run, and
    a single probability~$p\in\{1,0.8,0.6,0.4,0.2,0.1,0.05,0.01\}$ was fixed for each run 
    to simulate total asynchrony.
    At each timestep $k\in\mathbb{N}$, 
    processor~$n \in [N]$ draws three random values $p^n_{\text{up}},p^n_\text{comm},p^n_\text{recv}\in[0,1]$ from a uniform distribution.
    If $p^n_\text{up} < p$, then processor $n$ updates its own decision variable $W^n_n$ at time step $k$. 
    If $p_\text{comm}^n < p$, then processor $n$ sends its most recent copy of $W^n_n$ to all other processors $j\in[N]\backslash \{n\}$ at timestep $k$. 
    If $p_\text{recv}^n < p$, 
    then processor~$n$ receives the communications from all other processors that have been sent but not yet received. 
    It uses these received values to overwrite previous values of those processors' decision variables that it had stored onboard.
    
    In these simulations, each processor is initialized with a random matrix of size $K\times D$ drawn from a scaled normal distribution. 
    We compare Algorithm \ref{alg:TAGM} to totally asynchronous versions of gradient
    descent, the heavy ball algorithm, and Nesterov's accelerated gradient algorithm. 
    The four algorithms were run concurrently, 
    so that each algorithm has the same initial state and each processor $n\in[N]$ updates, communicates, and receives information with the exact same timing within each algorithm.

    The training procedure described above yielded [$76\%$-$80\%$] accuracy on the test set for all algorithms.
    Since the objective~$f$ in \eqref{eq:newsgroupF} satisfies Assumption \ref{ass:H}, 
    we know that~$f$ is strongly convex and hence it has a unique minimizer. 
    Each algorithm converges to that same minimizer, and 
    we therefore compare convergence rates. 
    For all algorithms, the step size is $\gamma=1.0$.
    The parameters $\beta$ and $\lambda$ are offset by $-0.001$ from the upper bounds needed for the guaranteed convergence of each algorithm (these bounds can be found in \cite{bertsekas1989parallel} for GD, \cite{hustig2023totally} for HB, and \cite{pond2024} for NAG.) 
    The values of $(\gamma,\lambda,\beta)$ are shown in Table~\ref{tab:param}.

    \begin{table}[ht]
        \centering
        \begin{tabular}{c| cccc}
            \toprule
            \textbf{Parameter} & \textbf{GD} & \textbf{HB} & \textbf{NAG} & \textbf{GM} \\
            \midrule
             $\gamma$ & 1.000 & 1.000 & 1.000 & 1.000 \\
             $\beta$ & -- & 0.249 & 0.499 & 0.248 \\
             $\lambda$ & -- & -- & 0.499 & 0.747 \\
             \bottomrule
        \end{tabular}
        \caption{Parameter values for each totally asynchronous algorithm.}
        \label{tab:param}
    \end{table}
    
    \begin{table}[ht]
        \centering
        \begin{tabular}{c|cccc}
            \toprule
            \multirow{2}{*}{\textbf{Probability}} &
            \multicolumn{4}{c}{\textbf{Number of Iterations}} 
            \\ 
            & \textbf{GD} & \textbf{HB} & \textbf{NAG} &\textbf{GM}  \\
            \midrule
            1.0 & 532 & 468 & 403 & 335 \\
            0.8 & 666 & 585 & 504 & 420 \\
            0.6 & 894 & 787 & 679 & 564 \\
            0.4 & 1343 & 1184 & 1022 & 850 \\
            0.2 & 2702 & 2387 & 2058 & 1726 \\
            0.1 & 5396 & 4765 & 4119  & 3429 \\
            0.05 & 10,758 & 9460 & 8150 & 6766 \\
            0.01 & 53,346 & 47,043 & 40,783 & 33,903 \\
            \bottomrule
        \end{tabular}
        \caption{Number of iterations needed to produce an iterate 
        within $\epsilon=10^{-3}$ of the optimum.}
        \label{tab:results}
    \end{table}

\begin{figure}[H]
    \centering
        \includegraphics[trim=.5cm 0cm 1.2cm .5cm, clip, width=\linewidth]{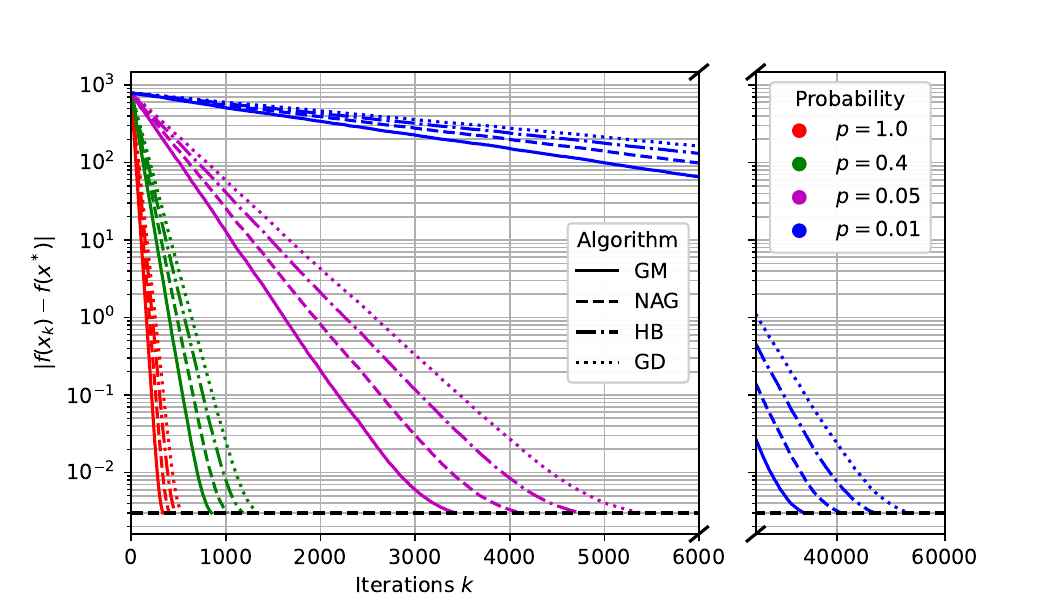}
        \caption{Cost convergence of totally asynchronous algorithms: GM (solid), NAG (dashed), HB (dash-dotted), and GD (dotted).}
        \label{fig:sims}
    \end{figure}

    Figure~\ref{fig:sims} shows the convergence of the cost
    for each algorithm for all values of $p$. The cost at each
    iteration was evaluated at the point~$(W^1_1(k), W^2_2(k), \ldots, W^N_N(k))$. 
    The numbers of iterations required 
    for each algorithm to produce such a point within
    a~$2$-norm distance
    $\epsilon = 10^{-3}$ of an optimum
    are shown in Table \ref{tab:results}.
    Figure~\ref{fig:sims} and Table~\ref{tab:results}
    demonstrate the improved performance of Algorithm 1 compared to each other 
    totally asynchronous algorithm. 
    Across all values of~$p$, 
    the decrease in required iterations from NAG to GM was $16\%$-$17\%$, 
    the decrease
    from HB to GM was $28\%$-$29\%$, and 
    the decrease
    from GD to GM was~$36\%$-$37\%$. 
    These results indicate that a
    significant reduction in computation time is attained by using the proposed GM algorithm,
    whether delays are short or long.

\section{Conclusion}

We presented a totally asynchronous generalized momentum algorithm
that subsumes several prior algorithms. 
We proved that it converges linearly in the number of operation cycles completed and 
we characterized the number of computations and communications that must be performed 
for processors' iterates to get within a desired error ball about an optimum. 
Simulations verified that this algorithm significantly outperforms existing ones. 
Future work will seek to develop rules to change algorithm parameters 
online to accelerate convergence even further.

\appendix

\subsection{Proof of Lemma~\ref{lem:sssFP}} \label{app:lem1proof}
The function $f$ has a unique minimizer $x^\star=\col{x_1^\star,\dots}{x_n^\star}\in\X$ defined in~\eqref{eq:xstardef}.
    By Proposition 5.7 in \cite[Section 3.5.5]{bertsekas1989parallel}, $x^\star$ 
    is the unique minimizer if and only if it 
    satisfies the variational inequality  
    \begin{equation}\label{eq:lem1_1}
        \langle x^i_i-x^\star_i,\nabla_if(x^\star)\rangle\geq 0
    \end{equation}
    for all $x^i_i \in \X_i$ and $i \in \V$.    
    The inequality in \eqref{eq:lem1_1} can be scaled by $\gamma > 0$ to show that 
    $x^{\star}$ is the unique minimizer if and only if it satisfies
    $\langle x^i_i-x^\star_i,\gamma\nabla_if(x^\star)\rangle\geq 0$.
    The term $\gamma\nabla_if(x^\star)$ is equivalently expressed as
    \begin{multline}\label{eq:lem1_2}
        \gamma\nabla_if(x^\star) = -x^\star_i + \gamma \nabla_if(x^\star+\lambda(x^\star-x^\star))\\-\beta(x^\star_i-x^\star_i)+x^\star_i,
    \end{multline}
    where we have added zero on the right-hand side.
    Substituting \eqref{eq:lem1_2} into \eqref{eq:lem1_1} and multiplying by~$-1$ yields
    \begin{multline}\label{eq:lem1_3}
        \langle x^i_i-x^\star_i,x^\star_i-\gamma\nabla_if(x^\star+\lambda(x^\star-x^\star)) \\ +\beta(x^\star_i-x^\star_i)-x^\star_i)\rangle \leq0
    \end{multline}
    for all $x^i_i\in\X_i$ and $i\in\V$.
    By the projection theorem in \cite{bertsekas1989parallel}, for a nonempty, closed, and convex set $\X$ and a continuously differentiable and convex function $f$, a vector $w\in\X$ satisfies $w=\Pi_{\X}[v]$ with $v\in\R^n$ if and only if $\langle y-w,v-w\rangle\leq0$ for all $y\in\X$.
    In \eqref{eq:lem1_3}, let $w=x^\star_i$, $v=x^\star_i-\gamma\nabla_if(x^\star+\lambda(x^\star-x^\star))+\beta(x^\star_i-x^\star_i)$, and $y=x^i_i$.
    Then
    \begin{equation}\label{eq:lem1_4}
        x^\star_i = \Pi_{\X_i}[x^\star_i-\gamma\nabla_i f(x^\star+\lambda(x^\star-x^\star))+\beta(x^\star_i-x^\star_i)],
    \end{equation}
    which shows~that $x^\star_i=\tilde{u}^i_x(x^\star,x^\star)$ for all $i\in\V$.
    Then~$(x^\star,x^\star)$ is a fixed point of the update law \eqref{eq:sssuxi} for all $i\in\V$.

    Consider now the update law~\eqref{eq:sssuyi}, where $y_i^i(l+1)=\tilde{u}^i_y(x^i(l),y^i(l))= x^i_i(l)$. 
    Plugging in $(x^\star,x^\star)$ gives $x^\star_i=\tilde{u}^i_y(x^\star,x^\star)$ for all $i\in\V$.
    Therefore, $(x^\star,x^\star)$ is a fixed point of the update law \eqref{eq:sssuyi} for all $i\in\V$.

    Thus, $z^\star_i=(x^\star_i,x_i^\star)\in\Z_i$ is a fixed point of the update law in~\eqref{eq:sssui} for all $i\in\V$ because $z^\star_i=\tilde{u}^i(z^\star)$ for all $i\in\V$.

    \subsection{Technical Lemmas} \label{app:tech}
    Toward proving Theorem~\ref{thrm:sssConv} in Appendix~\ref{app:main}, we first present two
    technical lemmas. 

    \begin{lemma} \label{lem:new}
     Consider the problem in~\eqref{eq:prob1} 
     and let Assumptions~\ref{ass:X}-\ref{ass:H} hold. Define~$a^i=(a^i_1,a^i_2)$ and $c^i=(c^i_1,c^i_2)$.
     If~$(\gamma,\lambda,\beta)\in C_1\cup C_2$, then 
    \begin{equation} \label{eq:lem2first}
        1- \gamma(1+\lambda)H_{ii}\big(b^i_1+\lambda(b^i_1-b^i_2)\big) + \beta > 0
    \end{equation}
    and
    \begin{equation} \label{eq:lem2second}
        \gamma\lambda H_{ii}\big(b^i_1+\lambda(b^i_1-b^i_2)\big) - \beta < 0
    \end{equation}
    for some $b^i_1,b^i_2 \in \mathbb{R}^n$ such that $b^i_{1,j} =\rho a^i_{1,j}+(1-\rho)c^i_{1,j}$ and~$b^i_{2,j}=\rho a^i_{2,j}+(1-\rho)c^i_{2,j}$ for all~$j \in \V_i$, where~$\rho \in [0, 1]$.
\end{lemma}
    \begin{proof}
        Recall the multi-variate Mean Value Theorem (MVT)~\cite[Theorem 1.9.1]{hubbard2002vector}:
    For a continuously differentiable function $h:\mathbb{R}^n\rightarrow\mathbb{R}$, there exists some $\rho\in[0,1]$ such that
    \begin{equation}
        h(v)-h(w)=\big\langle \nabla h\big(\rho w+(1-\rho)v\big),v-w\big\rangle
    \end{equation}
    for~$v, w \in \mathbb{R}^n$. 
    We will apply the MVT to the mapping~$\hat{u}^i_{\tr}$ from~\eqref{eq:utrue}.
    We use~$a^i=(a^i_1,a^i_2)$ and $c^i=(c^i_1,c^i_2)$ to find
    \begin{multline}\label{eq:sssproof1_lem}
        \hat u^i_\tr(c^i)-\hat u^i_\tr(a^i) = \sum_{j=1}^n \frac{\partial \hat u^i_\tr(b^i)}{\partial x_j}(c^i_{1,j}-a^i_{1,j})
        \\ + \sum_{j=1}^n \frac{\partial \hat u^i_\tr(b^i)}{\partial y_j}(c^i_{2,j}-a^i_{2,j})
    \end{multline}
    where~$b^i_1,b^i_2 \in \mathbb{R}^n$ satisfy
    $b^i_{1,j} =\rho a^i_{1,j}+(1-\rho)c^i_{1,j}$ and~$b^i_{2,j}=\rho a^i_{2,j}+(1-\rho)c^i_{2,j}$ for all~$j \in \V$
    and~$\rho \in [0, 1]$. 
    For~$j \neq i$, the partial derivatives in~\eqref{eq:sssproof1_lem} are
    \begin{align}
        \frac{\partial \hat u^i_\tr(b^i)}{\partial x_i}&= 1-\gamma(1+\lambda)\nabla^2_if\big(b^i_1+\lambda(b^i_1-b^i_2)\big)+\beta
        \\ \frac{\partial \hat u^i_\tr(b^i)}{\partial x_j} &= -\gamma(1+\lambda)\nabla_j\nabla_if\big(b^i_1+\lambda(b^i_1-b^i_2)\big)
        \\ \frac{\partial \hat u^i_\tr(b^i)}{\partial y_i} &= \gamma\lambda \nabla^2_i f \big(b^i_1 + \lambda (b^i_1-b^i_2)\big)- \beta
        \\ \frac{\partial \hat u^i_\tr(b^i)}{\partial y_j} &= \gamma\lambda \nabla_j\nabla_i f \big( b^i_1 + \lambda(b^i_1-b^i_2)\big). 
    \end{align}    
    Substituting these partial derivatives into \eqref{eq:sssproof1_lem} yields
    \begin{multline}
        \hat u^i_\tr(c^i)-\hat u^i_\tr(a^i)\\ = \Big(1 - \gamma(1+\lambda)\nabla_i^2f\big( b^i_1 + \lambda(b^i_1-b^i_2)\big)
        +\beta\Big)(c^i_{1,i}-a^i_{1,i}) 
        \\+\sum_{\substack{j=1\\j\not=i}}^n\Big({-\gamma}(1+\lambda) \nabla_j\nabla_i f \big(b^i_1 + \lambda(b^i_1-b^i_2)\big)\Big)(c^i_{1,j}-a^i_{1,j})
        \\ + \Big(\gamma\lambda\nabla^2_if\big(b^i_1+\lambda(b^i_1-b^i_2)\big)-\beta\Big)(c^i_{2,i}-a^i_{2,i})
        \\ + \sum_{\substack{j=1\\j\not=i}}^n\gamma\lambda \nabla_j\nabla_if\big(b^i_1+\lambda(b^i_1-b^i_2)\big)(c^i_{2,j}-a^i_{2,j}).
    \end{multline}
    We can replace $\nabla_i^2f$ with $H_{ii}$ and $\nabla_j\nabla_if$ with $H_{ij}$, 
    where~$H(x)=\nabla^2f(x)$ is the Hessian of the objective function~$f$. 
    Then we take the absolute value of both sides and apply the triangle inequality to obtain
    \begin{multline}\label{eq:sssproof2_lem}
        \vert \hat u^i_\tr(c^i)-\hat u^i_\tr(a^i)\vert
        \\ \leq \big\vert 1-\gamma(1+\lambda)H_{ii}\big(b^i_1+\lambda(b^i_1-b^i_2)\big)+\beta\big\vert\vert c^i_{1,i}-a^i_{1,i}\vert
        \\ +\sum_{\substack{j=1\\j\not=i}}^n\big\vert {-\gamma}(1+\lambda)H_{ij}\big(b^i_1+\lambda(b^i_1-b^i_2)\big)\big\vert\vert c^i_{1,j}-a^i_{1,j}\vert
        \\+\big\vert \gamma\lambda H_{ii}\big(b^i_1+\lambda(b^i_1-b^i_2)\big)-\beta \big\vert\vert c^i_{2,i}-a^i_{2,i}\vert
        \\ +\sum_{\substack{j=1\\j\not=i}}^n\big\vert \gamma\lambda H_{ij} \big(b^i_1+\lambda(b^i_1-b^i_2)\big)\big\vert\vert c^i_{2,j}-a^i_{2,j}\vert.
    \end{multline}
    
    Since $(\gamma, \lambda, \beta) \in C_1\cup C_2$, we consider two cases.
    First, if~$(\gamma,\lambda,\beta)\in C_1$, then
    \begin{align}
        1 - \gamma(1 &+ \lambda)H_{ii}(b^i_1 +\lambda(b^i_1-b^i_2)) + \beta 
    \\&> 1-\frac{\beta(1+\lambda)}{\lambda \max\limits_{i\in\V}\max\limits_{\eta\in\X}\vert H_{ii}(\eta)\vert}\max\limits_{i\in\V}\max\limits_{\eta\in\X}\vert H_{ii}(\eta)\vert +\beta
    \\&= 1-\frac{\beta(1+\lambda)}{\lambda}+\beta
    \\&= 1-\frac{\beta}{\lambda}\geq 0,
    \end{align}
    where the first inequality is obtained by replacing $\gamma$ with its upper bound $\frac{\beta}{\lambda \max\limits_{i\in\V}\max\limits_{\eta\in\X}\vert H_{ii}(\eta)\vert}$ and replacing $H_{ii}\big(b^i_1+\lambda(b^i_1-b^i_2)\big)$ with its upper bound $\max\limits_{i\in\V}\max\limits_{\eta\in\X}\vert H_{ii}(\eta)\vert$. The last inequality is obtained by using $\frac{\beta}{\lambda}\leq1$. We also find that 
    \begin{align}
        \gamma\lambda H_{ii}&\big(b^i_1+\lambda(b^i_1-b^i_2)\big)-\beta
    \\&< \frac{\beta\lambda}{\lambda \max\limits_{i\in\V}\max\limits_{\eta\in\X}\vert H_{ii}(\eta)\vert}\max\limits_{i\in\V}\max\limits_{\eta\in\X}\vert H_{ii}(\eta)\vert -\beta
    \\&= \beta - \beta = 0, 
    \end{align}
    where the first inequality is obtained by replacing $\gamma$ with its upper bound $\frac{\beta}{\lambda \max\limits_{i\in\V}\max\limits_{\eta\in\X}\vert H_{ii}(\eta)\vert}$ and replacing $H_{ii}\big(b^i_1+\lambda(b^i_1-b^i_2)\big)$ with its upper bound $\max\limits_{i\in\V}\max\limits_{\eta\in\X}\vert H_{ii}(\eta)\vert$.
    Then~\eqref{eq:lem2first} and~\eqref{eq:lem2second} hold for~$(\gamma,\lambda,\beta)\in C_1$. 
    
    Now we consider~$(\gamma,\lambda,\beta)\in C_2$. Then 
    \begin{align}
        1-\gamma(&1+\lambda)H_{ii}(b^i_1 +\lambda(b^i_1-b^i_2)) + \beta 
    \\&> 1-\frac{1+\lambda}{\max\limits_{i\in\V}\max\limits_{\eta\in\X}\vert H_{ii}(\eta)\vert}\max\limits_{i\in\V}\max\limits_{\eta\in\X}\vert H_{ii}(\eta)\vert +\beta
    \\&= 1-(1+\lambda)+\beta
    \\&= \beta-\lambda\geq 0,
    \end{align}
    where the first inequality is obtained by replacing $\gamma$ with its upper bound $\frac{1}{ \max\limits_{i\in\V}\max\limits_{\eta\in\X}\vert H_{ii}(\eta)\vert}$ and replacing $H_{ii}\big(b^i_1+\lambda(b^i_1-b^i_2)\big)$ with its upper bound $\max\limits_{i\in\V}\max\limits_{\eta\in\X}\vert H_{ii}(\eta)\vert$. The last inequality is obtained by using $\beta \geq \lambda$. We also find 
    \begin{align}
        \gamma\lambda H_{ii}\big(&b^i_1+\lambda(b^i_1-b^i_2)\big)-\beta
    \\&< \frac{\lambda}{ \max\limits_{i\in\V}\max\limits_{\eta\in\X}\vert H_{ii}(\eta)\vert}\max\limits_{i\in\V}\max\limits_{\eta\in\X}\vert H_{ii}(\eta)\vert -\beta
    \\&= \lambda - \beta \leq 0 
    \end{align}
    where the first inequality is obtained by replacing $\gamma$ with its upper bound $\frac{1}{ \max\limits_{i\in\V}\max\limits_{\eta\in\X}\vert H_{ii}(\eta)\vert}$ and replacing $H_{ii}\big(b^i_1+\lambda(b^i_1-b^i_2)\big)$ with its upper bound $\max\limits_{i\in\V}\max\limits_{\eta\in\X}\vert H_{ii}(\eta)\vert$. The last inequality is obtained by using $\lambda\leq\beta$.
    Then~\eqref{eq:lem2first} and~\eqref{eq:lem2second} hold for~$(\gamma,\lambda,\beta)\in C_2$. 
    \end{proof}
    
    \begin{lemma}\label{lem:alphabound}
        Consider the problem in~\eqref{eq:prob1} and let Assumptions \ref{ass:X}-\ref{ass:H} hold.
        Define
        \begin{align}\label{eq:alpha1}
        \alpha_1&= \big(1+\beta-\gamma\mu(1+\lambda)\big)^2 
        \\ &\quad\quad\quad + (\beta-\gamma\lambda \mu)\big(2+\beta-\gamma\mu(1+\lambda)\big)
    \end{align}
    and
    \begin{align}\label{eq:alpha2}
        \alpha_2 = 1-\gamma\mu+2(\beta-\lambda\gamma\mu).
    \end{align}
        If $(\gamma,\lambda,\beta)\in C_1\cup C_2$, then $\alpha=\max\{\alpha_1,\alpha_2\}$ 
        satisfies $\alpha\in(0,1)$.
    \end{lemma}
    \begin{proof}
        Since $(\gamma,\lambda,\beta)\in C_1\cup C_2$, we consider two cases.
    First, if $(\gamma,\lambda,\beta)\in C_1$, then
    \begin{align}
        \alpha_2 &= 1 + 2\beta - \gamma\mu(1+2\lambda)
        \\ &\leq 1 - \gamma\mu + 2\lambda(1-\gamma\mu)
        \\ &< 1 - \gamma\mu + \frac{\gamma\mu}{1-\gamma\mu}(1-\gamma\mu)
        \\ &= 1 -\gamma\mu + \gamma\mu =1,
    \end{align}
    where the first inequality is obtained by replacing $\beta$ with its upper bound $\lambda$ and the strict inequality is obtained from replacing $\lambda$ with its upper bound $\frac{\gamma\mu}{2(1-\gamma\mu)}$.
    Then
    \begin{align}
        \alpha_1 &= \big(1+\beta-\gamma\mu(1+\lambda)\big)^2
        \\ &\quad\quad\quad+(\beta-\gamma\lambda\mu)\big(2+\beta-\gamma\mu(1+\lambda)\big)
        \\&\leq \big(1 -\gamma\mu+\lambda(1-\gamma\mu)\big)^2
        \\ &\quad\quad\quad + \lambda(1-\gamma\mu)\big(2-\gamma\mu + \lambda(1-\gamma\mu)\big),
        \\ &<\Big(1-\gamma\mu + \frac{\gamma\mu(1-\gamma\mu)}{2(1-\gamma\mu)}\Big)^2
        \\ &\quad\quad\quad + \frac{\gamma\mu(1-\gamma\mu)}{2(1-\gamma\mu)}\Big(2 - \gamma\mu + \frac{\gamma\mu(1-\gamma\mu)}{2(1-\gamma\mu)}\Big)
        \\ &= \big(1-\gamma\mu + \frac{\gamma\mu}{2}\big)^2 + \frac{\gamma\mu}{2}\big(2-\gamma\mu +\frac{\gamma\mu}{2}\big)
        \\ &= \big(1-\frac{\gamma\mu}{2}\big)^2+ \frac{\gamma\mu}{2}\big(2-\frac{\gamma\mu}{2}\big)
        \\ &= 1 + \frac{(\gamma\mu)^2}{4}- \gamma\mu + \gamma\mu - \frac{(\gamma\mu)^2}{4} = 1,
    \end{align}
    where the first inequality is obtained by replacing $\beta$ with its upper bound $\lambda$ and the strict inequality comes from replacing $\lambda$ with its upper bound $\frac{\gamma\mu}{2(1-\gamma\mu)}$.
        
    Now consider the lower bound of $\alpha_2$ if $(\gamma,\lambda,\beta)\in C_1$, namely
    \begin{align}
        \alpha_2 &= 1 + 2\beta-\gamma\mu(1+2\lambda)
        \\&= 1 -\gamma\mu + 2(\beta-\gamma\lambda\mu)
        \\ &\geq 1-\gamma\mu,
    \end{align}
    where the inequality comes from multiplying the upper bound $\gamma< \frac{\beta}{\lambda\max\limits_{i\in\V}\max\limits_{\eta\in\X}\vert H_{ii}(\eta)\vert}$ by $\mu$ 
    to get $\gamma\mu< \frac{\beta}{\lambda}$ and then rearranging to get $\beta-\gamma\lambda\mu>0$.
    Using~$\beta \leq \lambda$ we see that~$\gamma\mu < \frac{\beta}{\lambda} \leq 1$ and 
    \begin{align}
        \alpha_2 &= 1-\gamma\mu>0.
    \end{align}
    Next, if $(\gamma,\lambda,\beta)\in C_1$, then we find the lower bound
    \begin{align}
        \alpha_1 &=\big(1 + \beta -\gamma\mu(1+\lambda)\big)^2 
        \\ &\quad\quad\quad + (\beta - \gamma\lambda\mu)\big(2+\beta-\gamma\mu(1+\lambda)\big)
        \\ &=\big((1-\gamma\mu) + (\beta-\gamma\lambda\mu)\big)^2
        \\ &\quad\quad\quad + (\beta-\gamma\lambda\mu)\big((2 -\gamma\mu) + (\beta-\gamma\lambda\mu)\big)^2
        >0,
    \end{align}
    where the inequalities~$1 - \gamma\mu > 0$, $2 - \gamma\mu > 0$, and
    $\beta - \gamma\lambda\mu > 0$ imply that the whole expression is strictly positive. 
    Then $\alpha=\max\{\alpha_1,\alpha_2\}\in (0,1)$ holds for $(\gamma,\lambda,\beta)\in C_1$.
      
    Next, if $(\gamma,\lambda,\beta)\in C_2$, then
    \begin{align}
        \alpha_2 &= 1+2\beta-\gamma\mu(1+2\lambda)
        \\ &< 1 +\gamma\mu(1+2\lambda) - \gamma\mu(1+2\lambda)=1,
    \end{align}
    where we have replaced $\beta$ with its strict upper bound $\frac{1}{2}\gamma\mu(1+2\lambda)$.
    For $\alpha_1$, we have the upper bound
    \begin{align}
        \alpha_1 &= \big(1+\beta-\gamma\mu(1+\lambda)\big)^2 
        \\ &\quad\quad\quad +  (\beta-\gamma\lambda\mu)\big(2+\beta-\gamma\mu(1+\lambda)\big)
        \\ &< \big(1 +\frac{1}{2}\gamma\mu+\lambda \gamma\mu -\gamma\mu-\lambda\gamma\mu\big)^2
        \\ &\quad\quad\quad + \big( \frac{1}{2}\gamma\mu+\lambda\gamma\mu-\lambda\gamma\mu\big)
        \\ & \quad\quad\quad \cdot\big(2 + \frac{1}{2}\gamma\mu+\lambda \gamma\mu - \gamma\mu-\lambda\gamma\mu\big)
        \\ &= \big( 1- \frac{1}{2}\gamma\mu\big)^2+ \frac{1}{2}\gamma\mu\big(2-\frac{1}{2}\gamma\mu\big)
        \\ &= 1 + \frac{1}{4}(\gamma\mu)^2 -\gamma\mu+\gamma\mu - \frac{1}{4}(\gamma\mu)^2=1,
    \end{align}
    where the inequality comes from the bound~$\beta<\frac{1}{2}\gamma\mu(1+2\lambda)$.
    We lower bound $\alpha_2$ as    
    \begin{align}
        \alpha_2 &= 1 + 2\beta-\gamma\mu(1+2\lambda)
            \\ &> 1+ 2\beta-(1+2\lambda)
            \\ &= 2(\beta-\lambda)\geq0,
    \end{align}
    where the first inequality is found by replacing~$\gamma\mu$ with its upper bound $1$ and 
    the second inequality is due to $\beta\geq\lambda$.
    For $\alpha_1$ we have
    \begin{align}
        \alpha_1 &=\big(1+\beta-\gamma\mu(1+\lambda)\big)^2
        \\ &\quad\quad\quad + (\beta-\gamma\lambda\mu)\big(2+\beta-\gamma\mu(1+\lambda)\big)
        \\ &\geq \big(1 -\gamma\mu + \lambda(1-\gamma\mu)\big)^2
        \\ &\quad\quad\quad + \lambda(1-\gamma\mu)\big(2-\gamma\mu+\lambda(1-\gamma\mu)\big)
        \\ &\geq (1-\gamma\mu)^2 >0,
    \end{align}
    where the first inequality is found by replacing $\beta$ with its lower bound $\lambda$. 
    The second inequality is from replacing $\lambda$ with its lower bound $0$.
    The strict inequality is due to $1-\gamma\mu>0$.
    Then $\alpha=\max\{\alpha_1,\alpha_2\}\in(0,1)$ holds for $(\gamma,\lambda,\beta)\in C_2$
    as well. 
    \end{proof}

    \subsection{Proof of Theorem~\ref{thrm:sssConv}} \label{app:main}
    Lemma~\ref{lem:sssFP} showed $z^\star=(x^\star,x^\star)$ is the fixed point of~\eqref{eq:sssui} for all $i\in\V$
    and thus~$x^\star_i=\tilde{u}^i_x(z^\star)$ and $x^\star_i=\tilde{u}^i_y(z^\star)$.
    By definition of the $\infty-$norm, we have
    \begin{equation}
        \Vert x^\tr(l+1)-x^\star\Vert_\infty = \max_{i\in\V}\vert x^i_i(l+1)-x^\star_i\vert.
    \end{equation}
    Substituting the update for $x^i_i$ given in \eqref{eq:sssuxi} and the fixed-point property of $x^\star$, we have    
    \begin{multline}\label{eq:refback}
        \Vert x^\tr(l+1)-x^\star\Vert_\infty = \max_{i\in\V}\Big\vert \Pi_{\X_i}\big[x^i_i(l)
        \\-\gamma\nabla_if\Big(x^\tr(l)+\lambda\big(x^\tr(l)-y^\tr(l)\big)\Big)+\beta\big(x^i_i(l)-y^i_i(l)\big)\big]
        \\ - \Pi_{\X_i}\big[x^\star_i-\gamma\nabla_if\big(x^\star+\lambda(x^\star-x^\star)\big)+\beta(x^\star_i-x^\star_i)\big]\Big\vert,
    \end{multline}
    where we have replaced~$x^i(l)$ with $x^\tr(l)$ and $y^i(l)$ with $y^\tr(l)$, which is justified
    as follows. 
    Processor~$i$ does not communicate with processor~$m$ if~$m \not \in \V_i$. 
    This setup is unproblematic because the values of~$z^i_m$
    for~$m \not\in \V_i$ do not affect processor~$i$'s computations and
    it was noted below Definition~\ref{def:zstar} that the values of~$z^i_m$
    can therefore be set to arbitrary values. 
    We can set~$z^i_m(l) = z^{\tr}_m(l)$ because doing
    so does not change the results of processor~$i$'s computations and
    we choose to do so for ease of analysis. 
    
    The orthogonal projection is continuous and nonexpansive, i.e., $\vert \Pi_{\X_i}[v_2]-\Pi_{\X_i}[v_1]\vert\leq \vert v_2-v_1\vert$ for all $v_1,v_2\in\R$ \cite{bertsekas1989parallel}.
    Therefore, 
    \begin{multline}
        \Vert x^\tr(l+1)-x^\star\Vert_\infty \leq \max_{i\in\V}\Big\vert x^i_i(l)-\gamma\nabla_i f\Big(x^\tr(l)
        \\+\lambda\big(x^\tr(l)-y^\tr(l)\big)\Big) + \beta\big(x^i_i(l)-y^i_i(l)\big) \\- \big(x^\star_i-\gamma\nabla_if\big(x^\star + \lambda(x^\star-x^\star)\big)+\beta(x^\star_i-x^\star_i)\big)\Big\vert.
    \end{multline}
    By the definition of the true state of the network in \eqref{eq:utrue}, we have
    \begin{equation}\label{eq:sssproof12}
        \Vert x^\tr(l+1)-x^\star\Vert_\infty \leq \max_{i\in\V}\vert \hat{u}^i_\tr\big(z^\tr(l)\big)-\hat{u}^i_\tr(z^\star)\vert.
    \end{equation}

    Lemma~\ref{lem:new} shows that
    \begin{multline}\label{eq:sssproof3}
        \big\vert 1- \gamma(1+\lambda)H_{ii}\big(b^i_1+\lambda(b^i_1-b^i_2)\big) + \beta \big\vert \\= 
        1- \gamma(1+\lambda)H_{ii}\big(b^i_1+\lambda(b^i_1-b^i_2)\big) + \beta
    \end{multline}
    and
    \begin{multline}\label{eq:sssproof4}
        \big\vert \gamma\lambda H_{ii}\big(b^i_1+\lambda(b^i_1-b^i_2)\big)-\beta \big\vert \\= -\gamma\lambda H_{ii}\big(b^i_1+\lambda(b^i_1-b^i_2)\big)+\beta,
    \end{multline}
    for $(\gamma,\lambda,\beta)\in C_1\cup C_2$,
    which allows us to write
    \begin{multline}\label{eq:sssproof9}
        \vert \hat u^i_\tr(c^i)-\hat u^i_\tr(a^i)\vert 
        \\ \leq \big(1-\gamma(1+\lambda) H_{ii}\big(b^i_1 +\lambda (b^i_1-b^i_2)\big) + \beta\big)\vert c^i_{1,i}-a^i_{1,i}\vert
        \\+ \sum_{\substack{j=1\\j\not=i}}^n\big\vert {-\gamma}(1+\lambda)H_{ij}\big(b^i_1+\lambda(b^i_1-b^i_2)\big)\big\vert \vert c^i_{1,j}-a^i_{1,j}\vert
        \\ +\big({-\gamma}\lambda H_{ii}\big(b^i_1+\lambda(b^i_1-b^i_2)\big)+\beta\big)\vert c^i_{2,i}-a^i_{2,i}\vert
        \\ + \sum_{\substack{j=1\\j\not=i}}^n \big\vert \gamma\lambda H_{ii}\big(b^i_1+\lambda(b^i_1-b^i_2)\big)\big\vert \vert c^i_{2,j}-a^i_{2,j}\vert,
    \end{multline}
    where $a^i$, $b^i$, and $c^i$ are as defined in Lemma \ref{lem:new}.
    By the definition of the $\infty-$norm, $\Vert c^i_1-a^i_1\Vert_\infty =\max\limits_{j\in\V}\vert c^i_{1,j}-a^i_{1,j}\vert$ and $\Vert c^i_2-a^i_2\Vert_\infty = \max_{j\in\V}\vert c^i_{2,j}-a^i_{2,j}\vert$.
    Applying these equalities to \eqref{eq:sssproof9} gives
    \begin{multline}\label{eq:sssproof10}
        \vert \hat u^i_\tr(c^i) - \hat{u}^i_\tr(a^i)\vert
        \leq \Big(1- \gamma(1+\lambda) H_{ii}\big(b^i_1+\lambda(b^i_1-b^i_2)\big)
        \\+\beta + \gamma(1+\lambda) \sum_{\substack{j=1\\j\not=i}}^n\big\vert H_{ij}\big(b^i_1 + \lambda(b^i_1-b^i_2)\big)\big\vert\Big)\Vert c^i_1-a^i_1\Vert_\infty
        \\ + \Big({-\gamma}\lambda H_{ii}\big(b^i_1+\lambda(b^i_1-b^i_2)\big)+\beta \\+ \gamma\lambda\sum_{\substack{j=1\\j\not=i}}^{n}\big\vert H_{ij}\big(b^i_1
        +\lambda(b^i_1-b^i_2)\big)\big\vert\Big)\Vert c^i_2-a^i_2\Vert_\infty.
    \end{multline}
    Due to Assumption \ref{ass:H}, we have~$H_{ii}(\eta)\geq \mu + \sum_{j=1,j\not=i}^n\vert H_{ij}(\eta)\vert$ for all $\eta\in\X$ and $i\in\V$, where $\mu > 0$.
    Rearranging and negating this relationship yields
    \begin{equation}
        -H_{ii}(\eta)+\sum_{\substack{j=1\\j\not=i}}^n \vert H_{ij}(\eta)\vert \leq -\mu.
    \end{equation}
    Applying this inequality to \eqref{eq:sssproof10} yields 
    \begin{multline}\label{eq:sssproof11}
        \vert \hat u^i_\tr(c^i)-\hat u^i_\tr(a^i)\vert \leq \big( 1+ \beta -\gamma\mu(1+\lambda)\big)\Vert c^i_1-a^i_1\Vert_\infty
        \\ + (\beta - \gamma\lambda\mu)\Vert c^i_2-a^i_2\Vert_\infty.
    \end{multline}
    
    We now set $c^i=z^\tr(l)$ and $a^i=z^\star$, which implies that $c^i_1=x^\tr(l)$, $c^i_2=y^\tr(l)$, and $a^i_1=a^i_2=x^\star$.
    Then \eqref{eq:sssproof11} becomes
    \begin{multline}\label{eq:sssproof13}
        \vert \hat u^i_\tr(z^\tr(l))-\hat u^i_\tr(z^\star)\vert \\\leq \big( 1+ \beta -\gamma\mu(1+\lambda)\big)\Vert x^\tr(l)-x^\star\Vert_\infty
        \\+ (\beta - \gamma\lambda\mu)\Vert y^\tr(l)-x^\star\Vert_\infty,
    \end{multline}
    which holds for all $i\in\V$.
    We substitute \eqref{eq:sssproof13} into~\eqref{eq:sssproof12} to reach
    \begin{multline}\label{eq:sssproof14}
        \Vert x^\tr(l+1)-x^\star\Vert_\infty \leq \max_{i\in\V} \big( 1+ \beta -\gamma\mu(1+\lambda)\big)
        \\ \cdot\Vert x^\tr(l)-x^\star\Vert_\infty
        + (\beta - \gamma\lambda\mu)\Vert y^\tr(l)-x^\star\Vert_\infty
        \\ = \big( 1+ \beta -\gamma\mu(1+\lambda)\big)
         \Vert x^\tr(l)-x^\star\Vert_\infty
        \\+ (\beta - \gamma\lambda\mu)\Vert y^\tr(l)-x^\star\Vert_\infty,
    \end{multline}
    where equality holds because the argument of the maximum operator is independent of the index $i\in\V$.
    Since \eqref{eq:sssproof14} holds for all $l\in\mathbb{N}$, we also have 
    \begin{multline}\label{eq:sssproof15}
        \Vert x^\tr(l)-x^\star\Vert_\infty \leq \big( 1+ \beta -\gamma\mu(1+\lambda)\big)
         \Vert x^\tr(l-1)-x^\star\Vert_\infty
        \\+ (\beta - \gamma\lambda\mu)\Vert y^\tr(l-1)-x^\star\Vert_\infty.
    \end{multline}
    Substituting~$y^\tr(l) = x^\tr(l-1)$ in \eqref{eq:sssproof14} yields
    \begin{multline}\label{eq:sssproof16}
        \Vert x^\tr(l+1)-x^\star\Vert_\infty \leq \big(1+\beta-\gamma\mu(1+\lambda)\big)\Vert x^\tr(l)-x^\star\Vert_\infty
        \\ +(\beta -\gamma\lambda\mu)\Vert x^\tr(l-1)-x^\star\Vert_\infty.
    \end{multline}
    Now we substitute \eqref{eq:sssproof15} into \eqref{eq:sssproof16} to obtain
    \begin{multline}\label{eq:sssproof17}
        \Vert x^\tr(l+1)-x^\star\Vert_\infty 
        \leq \big(1+\beta -\gamma\mu(1+\lambda)\big)\\\cdot\Big(\big(1+\beta-\gamma\mu(1+\lambda)\big)\Vert x^\tr(l-1)-x^\star\Vert_\infty
        \\+(\beta-\gamma\lambda\mu)\Vert y^\tr(l-1)-x^\star\Vert_\infty\Big) \\+ (\beta-\gamma\lambda\mu)\Vert x^\tr(l-1)-x^\star\Vert_\infty
        \\ = \Big(\big(1+\beta-\gamma\mu(1+\lambda)\big)^2+\beta-\gamma\lambda\mu\Big)\Vert x^\tr(l-1)-x^\star\Vert_\infty
        \\ \big(1+\beta-\gamma\mu(1+\lambda)\big)(\beta-\gamma\lambda\mu)\Vert y^\tr(l-1)-x^\star\Vert_\infty.
    \end{multline}
    We also substitute $y^\tr(l+1)=x^\tr(l)$ into \eqref{eq:sssproof15} to find 
    \begin{multline}\label{eq:sssproof18}
        \Vert y^\tr(l+1)-x^\star\Vert_\infty\leq \big(1+\beta-\gamma\mu(1+\lambda)\big)
        \Vert x^\tr(l-1)-x^\star\Vert_\infty
         \\+ (\beta-\gamma\lambda\mu)\Vert y^\tr(l-1)-x^\star\Vert_\infty. 
    \end{multline}
    By the definition of the $\infty-$norm, for $z^\tr(l-1)=(x^\tr(l-1),y^\tr(l-1))$ we have
    \begin{align}
        \Vert x^\tr(l-1)-x^\star\Vert_\infty \leq \Vert z^\tr(l-1)-x^\star\Vert_\infty
        \\ \Vert y^\tr(l-1)-x^\star\Vert_\infty \leq \Vert z^\tr(l-1)-x^\star\Vert_\infty.
    \end{align}
    Therefore, \eqref{eq:sssproof17} can be further bounded as
    \begin{multline}\label{eq:sssproof19}
        \Vert x^\tr(l+1)-x^\star\Vert_\infty \leq \Big(\big(1+\beta-\gamma\mu(1+\lambda)\big)^2 \\+(\beta-\gamma\lambda\mu)\big(2+\beta-\gamma\mu(1+\lambda)\big)\Big)\Vert z^\tr(l-1)-x^\star\Vert_\infty
    \end{multline}
    and \eqref{eq:sssproof18} can be bounded as
    \begin{multline}\label{eq:sssproof20}
        \Vert y^\tr(l+1)-x^\star\Vert_\infty \leq\\ \big(1 + 2\beta -\gamma\mu(1+2\lambda) \big)\Vert z^\tr(l-1)-x^\star\Vert_\infty.
    \end{multline}
    Finally, we combine \eqref{eq:sssproof19} and \eqref{eq:sssproof20} by the definition of the $\infty-$norm to obtain 
    \begin{equation}
        \Vert z^\tr(l+1)-x^\star\Vert_\infty\leq \alpha\Vert z^\tr(l-1)-x^\star\Vert_\infty,
    \end{equation}
    for all $l\in\mathbb{N}$,
    where $\alpha \in (0, 1)$ is from
    Lemma~\ref{lem:alphabound}. 

    \subsection{Proof of Theorem~\ref{thrm:asymptotic}} \label{app:zsets}
    The asymptotic convergence of Algorithm~\ref{alg:TAGM} under total asynchrony follows from the satisfaction
    of the conditions in Lemma~\ref{lem:threepart}, which is what we prove. 
    
    (LLC) Lemma \ref{lem:alphabound} shows $\alpha\in(0,1)$. 
        Then~$\alpha^{k+1}< \alpha^k$ for all $k\in\mathbb N$ and
        $v \in \Z(k+1)$ implies that
        \begin{equation}
            \Vert v-z^\star\Vert_\infty\leq \alpha^{k+1}\Vert z(0)-z^\star\Vert_\infty  < \alpha^k\Vert z(0)-z^\star\Vert_\infty,
        \end{equation}
        and thus~$v \in \Z(k)$. 
        Then $\Z(k) \supset \Z(k+1)$ for all $k\in\mathbb N$.
        
        (SCC) 
        The mapping $h$ in~\eqref{eq:h} applies two iterations of the GM algorithm and 
        hence is an~$\alpha$-contraction with respect to the~$\infty$-norm by Theorem \ref{thrm:sssConv}. 
        A point~$z \in \Z(k)$ satisfies $\Vert z-z^\star\Vert\leq \alpha^k\Vert z(0)-z^\star\Vert_\infty$ and
        since $h$ is an $\infty$-norm contraction we have
        \begin{equation}
            \Vert h(z)-z^\star\Vert_\infty \leq \alpha\Vert z-z^\star\Vert_\infty \leq \alpha^{k+1}\Vert z(0)-z^\star\Vert_\infty.
        \end{equation}
        Then we have $h(z) \in \Z(k+1)$ for any~$k \in \mathbb{N}$ and $z \in \Z(k)$. 

        Consider a sequence $\{z_k\}_{k\in\mathbb N}$ with $z_k\in \Z(k)$ for all $k\in\mathbb{N}$.
        Since $\Z(k)$ satisfies the LLC, we know that $\Z(k-1)\supset\Z(k)$ for all $k\in\mathbb N$.
        By inspection of \eqref{eq:Zk}, each set $\Z(k)$ is closed.
        Therefore, from \cite{rockafellar2009variational}, we have $\lim\limits_{k\rightarrow\infty}\Z(k)=\bigcap_{k\in\mathbb N}\Z(k)$ and we find 
        \begin{equation}
            \lim_{k\rightarrow \infty}z_k\in \lim_{k\rightarrow \infty}\Z(k)= \bigcap_{k\in\mathbb N}\Z(k) = \{z^\star\}. 
        \end{equation}
        Then~$\lim\limits_{k\rightarrow\infty} z_k = z^\star$. 
        Here, $z^\star$ is a fixed point of each map $u^i$ by Lemma \ref{lem:dGMfp}. 
        From the definition of the map $h$ in \eqref{eq:h}, it is also a fixed point of~$h$. 
        
        (BCC) If $v\in \Z(k)$, then using the definition of the $\infty-$norm gives 
        \begin{equation}
            \Vert v- z^\star\Vert_\infty = \max_{i\in\V}\vert v_i-z^\star_i\vert \leq \alpha^k\Vert z(0)-z^\star\Vert_\infty
        \end{equation}
        and $\vert v_i-z^\star_i\vert \leq \alpha^k\Vert z(0)-z^\star\Vert_\infty$ for all $i\in \V$.
        Define
        \begin{equation}
            \Z_i(k) =\{v_i\in \Z_i : \vert v_i-z^\star_i\vert \leq \alpha^k\Vert z(0)-z^\star\Vert_\infty\}.
        \end{equation}
        Then~$v \in \Z(k)$ has $v_i \in \Z_i(k)$ for all $i \in \V$ and
        $\Z(k)=\Z_1(k)\times\dots\times \Z_n(k)$ holds for all $k\in\mathbb N$.    

        \subsection{Proof of Lemma~\ref{lem:invariant}} \label{app:invariant}
        Theorem~\ref{thrm:asymptotic}
        showed that the set $\Z(k)$ satisfies the SCC for all $k \in \mathbb N$.
    Formally, if $z \in \Z(k)$, then $h(z)\in\Z(k+1)$ for all $k\in\mathbb N$.
    From Assumption \ref{ass:TA}, there can be arbitrarily long delays between computing, sending, and receiving information, but there can be no permanent cessation of them.
    Therefore, there exists a time $k_1\in K^i$ when processor $i\in\V$ has 
    used~$z^i(k_1)\in\Z(k_1)$ to compute~$z^i_i(k_1+1)\in\Z_i(k_1+1)$.
    From Theorem \ref{thrm:asymptotic}, the sets $\Z(k)$ satisfy the LLC and thus~$\Z_i(k_1+1)\subset \Z_i(k_1)$.
    Therefore, $z^i_i(k_1+1) \in \Z_i(k_1 + 1) \subseteq \Z_i(k_1)$.
    Since~$z^i(k_1) \in \Z(k_1)$, we have~$z^i_j(k_1) \in \Z_j(k_1)$. 
    If no communications are received by processor~$i$ at time~$k_1$, then 
    each entry~$z^i_j(k_1+1) = z^i_j(k_1)$ for~$j \neq i$. 
    Then~$z^i_j(k_1+1) \in \Z_j(k_1)$ for~$j \neq i$
    and~$z^i_i(k_1+1) \in \Z_i(k_1)$ from above. 
    Then processor~$i$ has~$z^i(k_1+1) \in \Z(k_1)$ and~$\Z(k_1)$
    is invariant under processor~$i$'s computations for all~$i \in \V$.

    When communications do occur after time~$k_1$, they consist of processors sharing
    new values of decision variables that they have computed. These computations
    use elements of~$\Z(k_1)$ for the same reasoning as above and communicating them
    shares components of vectors that are also in~$\Z(k_1)$. 
    The set~$\Z(k_1)$ is therefore invariant under processor~$i$'s communications for all~$i \in \V$.  

    \subsection{Proof of Theorem~\ref{thrm:TAGM}} \label{app:linear}
    Let $k_i\in K^i$ be the first timestep that processor $i\in\V$ performs
    a computation. 
    Then~$\big(x^i_i(k_i+1),y^i_i(k_i+1)\big)\gets u^i\big(x^i(k_i),y^i(k_i)\big)$ for all $i \in \V$.
    By the BCC and the SCC in Lemma \ref{lem:threepart}, we have $\big(x^i_i(k_i+1),y^i_i(k_i+1)\big)\in \Z_i(1)$.
    By Assumption \ref{ass:TA}, there exists a time $k'=\max\limits_{i\in\V}k_i+1$ such that $\big(x^i_i(k'),y^i_i(k')\big)\in\Z_i(1)$ holds for all $i\in\V$.    
    Recall that $\big(x^i_i(k'),y^i_i(k')\big)\in \Z_i(0)$ still holds for all $i\in\V$ due to 
    the LLC in Lemma \ref{lem:threepart}.    
    Prior to processor $i\in\V$ receiving communications from any processors $j\in\V_i$, 
    the relation~$\big(x^i_j(k'),y^i_j(k')\big)\in\Z_j(0)$ holds.
    Similarly, processor~$i$ has~$\big(x^i_m(k'),y^i_m(k')\big)\in\Z_m(0)$
    for all $m \not\in \V_i$.

    Due to Assumption \ref{ass:TA}, there exists a timestep when 
    all processors will have sent messages to their neighbors and had those messages received by them. 
    Let $k_{i,j}\in R^i_j$ be the first time that processor $i\in\V$ receives $\big(x^j_j\big(\tau^i_j(k_{i,j})\big),y^j_j\big(\tau^i_j(k_{i,j})\big)\big)$ from processor~$j \in \V_i$.
    Then for~$k''=\max\limits_{i\in\V}\max\limits_{j\in\V_i} k_{i,j}$ 
    we have~$k'' \geq k'$ and
    by time~$k''$ each processor has received a communication from each of its neighbors. 
    For all processors outside of $i$'s neighborhood, that is, processors $m\not\in\V^i$, at every time step $(x^i_m(k),y^i_m(k))\gets (x^i_m(k-1),y^i_m(k-1)).$
    Therefore at time $k''$, we have $\big(x^i_j(k''),y^i_j(k'')\big)\in\Z_j(1)$ for all $i\in\V$ and $j\in\V_i$.
    Then processor~$i$'s local decision variable satisfies $\big(x^i(k''),y^i(k'')\big)\in\Z(1)$
    for all~$i \in \V$.
    Therefore, the first operation cycle has been completed and $\ops{k''}=1$.

    By the definition of $\Z(k)$ in \eqref{eq:Zk}, processor~$i$'s local 
    decision vector satisfies 
    \begin{equation}
        \Vert z^i(k'')-z^\star\Vert_\infty\leq \alpha\Vert z^i(0)-z^\star\Vert_\infty
    \end{equation}
    for all~$i \in \V$, where $\alpha^{\ops{k''}}=\alpha$ since $\ops{k''}=1$.
    By induction, we establish 
    \begin{equation}
        \Vert z^i(k)-z^\star\Vert_\infty\leq \alpha^{\ops{k}}\Vert z^i(0)-z^\star\Vert_\infty,
    \end{equation}
    for all $i\in\V$.
    Over all processors, we then have
    \begin{equation}
        \max_{i\in\V}\Vert z^i(k)-z^\star\Vert_\infty \leq \alpha^{\ops{k}}\max_{i\in\V}\Vert z^i(0)-z^\star\Vert_\infty,
    \end{equation}
    which is linear convergence of Algorithm \ref{alg:TAGM}.

    \subsection{Proof of Theorem~\ref{thrm:opcomplex}} \label{app:complex}
    Theorem~\ref{thrm:TAGM} gives
    \begin{equation}
        \max_{i\in\V}\Vert z^i(k)-z^\star\Vert_\infty \leq \alpha^{\ops{k}}\max_{i\in\V}\Vert z^i(0)-z^\star\Vert_\infty
    \end{equation}
    for all $k\in \mathbb N$.
    Then
    \begin{align}
        \max_{i\in\V}\Vert z^i(k)-z^\star\Vert_\infty &\leq \alpha^{\ops{k}}\max_{v_1,v_2\in\Z}\Vert v_1-v_2\Vert_\infty
        \\ &\leq \alpha^{\ops{k}}D.
    \end{align}
    We seek to enforce $\alpha^{\ops{k}}D \leq \epsilon$. 
    Solving for $\ops{k}$, we find
        $\ops{k}\log(\alpha) \leq \log (\epsilon/D)$.
    Since $\alpha \in (0,1)$, we have~$\log(\alpha) < 0$ and therefore
    \begin{align}
        \ops{k}\geq \frac{\log(\epsilon/D)}{\log(\alpha)} = \frac{\log(D/\epsilon)}{\log(1/\alpha)}.
    \end{align}
    Let $\rho = \frac{\log(D/\epsilon)}{\log(1/\alpha)}$.
    Since each processor must perform one computation and send~$\vert \V_i\vert$ communications per operation cycle, each processor must perform $\rho$ computations and send $\rho\vert \V_i\vert$ communications to get within 
    distance $\epsilon$ of the minimizer $z^\star$.

\bibliographystyle{ieeetr}
\bibliography{ref}

\end{document}